
\input amstex
\loadbold

\define\Pee{{\Bbb P}}
\define\Zee{{\Bbb Z}}
\define\Cee{{\Bbb C}}

\define\proof{\demo{Proof}}
\define\endproof{\qed\enddemo}

\define\theorem#1{\proclaim{Theorem #1}}
\define\lemma#1{\proclaim{Lemma #1}}
\define\proposition#1{\proclaim{Proposition #1}}
\define\corollary#1{\proclaim{Corollary #1}}
\define\claim#1{\proclaim{Claim #1}}

\define\section#1{\specialhead #1 \endspecialhead}
\define\ssection#1{\medskip\noindent{\bf #1}}

\documentstyle{amsppt}
\leftheadtext{}
\rightheadtext{}

\pageno=1
\topmatter
\title Quantum cohomology \\
of projective bundles over $\Pee^n$
\endtitle
\author {Zhenbo Qin$^1$ and Yongbin Ruan$^2$}
\endauthor
\address Department of Mathematics, Oklahoma State University,
Stillwater, OK 74078
\endaddress
\email  zq\@math.okstate.edu
\endemail
\address Department of Mathematics, University of Utah,
Salt Lake City, UT 84112
\endaddress
\email ruan\@math.utah.edu
\endemail
\thanks ${}^1$Partially supported by an NSF grant
\endthanks
\thanks ${}^2$Partially supported by an NSF grant and a Sloan fellowship
\endthanks
\endtopmatter

\TagsOnRight
\NoBlackBoxes
\document
\section{1. Introduction}

Quantum cohomology, proposed by Witten's study \cite{16} of
two dimensional nonlinear sigma models, plays a fundamental role
in understanding the phenomenon of mirror symmetry for Calabi-Yau manifolds.
This phenomenon was first observed by physicists motivated by
topological field theory. A topological field theory starts
with correlation functions. The correlation functions of sigma model
are linked with the intersection numbers of cycles in the moduli space
of holomorphic maps from Riemann surfaces to manifolds.
For some years, the mathematical construction of these correlation functions
remained to be a difficult problem because the moduli spaces
of holomorphic maps usually are not compact and may have wrong dimension.
The quantum cohomology theory was first
put on a firm mathematical footing by \cite{12,13} for semi-positive
symplectic manifolds (including  Fano and Calabi-Yau manifolds),
using the method of symplectic topology. Recently, an algebro-geometric
approach has been taken by \cite{8,9}. The results of \cite{12,13} have been
redone in the algebraic geometric setting for the case of homogeneous spaces.
The advantage of homogeneous spaces is that the moduli spaces of
holomorphic maps always have expected dimension and
their compactifications are nice. Beyond the homogeneous spaces,
one can not expect such nice properties for the moduli spaces.
The projective bundles are perhaps the simplest examples.
However, by developing sophisticated excessive intersection theory,
it is possible that the algebro-geometric method can work for
any projective manifolds.
In turn, it may shed new light to removing the semi-positive condition
in the symplectic setting.

Although we have a solid foundation for quantum cohomology theory at least for
semi-positive symplectic manifolds, the
calculation remains to be a difficult task.
So far, there are only a few examples
which have been computed, e.g., Grassmannian \cite{14},
some rational surfaces \cite{6},
flag varieties \cite{4}, some complete intersections \cite{3},
and the moduli space of stable bundles over Riemann surfaces \cite{15}.
One of the common feature for these examples is
that the relevant moduli spaces of rational curves have expected dimension.
Then, one can use the intersection theory. We should mention that there are
many predications based on mathematically unjustified mirror symmetry (for
Calabi-Yau 3-folds) and linear sigma model (for toric varieties). In this
paper, we attempt to determine the quantum cohomology of
projective bundles over the projective space $\Pee^n$. In contrast to the
previous examples, the relevant moduli spaces in our case frequently do not
have expected dimensions. It makes the calculation more difficult. We
overcome this difficulty by using  excessive intersection theory.

There are two main ingredients in our arguments. The first one is a result of
Siebert and Tian (the Theorem 2.2 in \cite{14}),
which says that if the ordinary cohomology $H^*(X; \Zee)$ of
a symplectic manifold $X$ with the symplectic form $\omega$
is the ring generated by $\alpha_1, \ldots, \alpha_s$
with the relations $f^1, \ldots, f^t$, then the quantum cohomology
$H^*_\omega(X; \Zee)$ of $X$ is the ring generated
by $\alpha_1, \ldots, \alpha_s$ with $t$ new relations
$f_\omega^1, \ldots, f_\omega^t$
where each new relation $f_\omega^i$ is just
the relation $f^i$ evaluated in the quantum cohomology ring structure.
It was known that the quantum product $\alpha \cdot \beta$ is
the deformation of ordinary cup product by the lower order terms
called quantum corrections. The second ingredient is that
under certain numerical conditions, most of the quantum corrections vanishes.
Moreover, the nontrivial quantum corrections seem to come from
Mori's extremal rays.

Let $V$ be a rank-$r$ bundle over $\Pee^n$,
and $\Pee(V)$ be the corresponding projective bundle.
Let $h$ and $\xi$ be the cohomology classes of a hyperplane in $\Pee^n$
and the tautological line bundle in $\Pee(V)$ respectively.
For simplicity, we make no distinction between $h$ and $\pi^*h$
where $\pi: \Pee(V) \to \Pee^n$ is the natural projection.
Denote the product of $i$ copies of $h$ and $j$ copies of $\xi$
in the ordinary cohomology ring by $h_i \xi_j$,
and the product of $i$ copies of $h$ and $j$ copies of $\xi$
in the quantum cohomology ring by $h^i \cdot \xi^j$.
For $i = 0, \ldots, r$, put $c_i(V) = c_i \cdot h_i$ for some integer $c_i$.
It is well known that $-K_{\Pee(V)} = (n + 1 - c_1)h + r\xi$
and the ordinary cohomology ring $H^*(\Pee(V); \Zee)$ is
the ring generated by $h$ and $\xi$ with the two relations:
$$h_{n+1} = 0 \qquad \text{and} \qquad \sum_{i=0}^r (-1)^i c_i \cdot h_i
\xi_{r-i} = 0. \eqno (1.1)$$
In particular, $H^{2(n+r-2)}(\Pee(V); \Zee)$ is generated by
$h_{n-1}\xi_{r-1}$ and $h_{n}\xi_{r-2}$,
and its Poincar\'e dual $H_2(\Pee(V); \Zee)$ is generated by
$(h_{n-1}\xi_{r-1})_*$ and $(h_{n}\xi_{r-2})_*$
where for $\alpha \in H^*(\Pee(V); \Zee)$,
$\alpha_*$ stands for its Poincar\'e dual. We have
$$-K_{\Pee(V)}(A) = a(n + 1 - c_1) + r \cdot \xi(A)
= a(n + 1 - c_1) + r(ac_1 + b) \eqno (1.2)$$
for $A = (a h_{n-1}\xi_{r-1}+ b h_{n}\xi_{r-2})_* \in H_2(\Pee(V); \Zee)$.

By definition, $V$ is an ample (respectively, nef) bundle if and only if
the tautological class $\xi$ is an ample (respectively, nef) divisor
on $\Pee(V)$. Assume that $V$ is ample such that either $c_1 \le (n+1)$
or $c_1 \le (n + r)$ and $V\otimes {\Cal O}_{\Pee^n}(-1)$ is nef.
Then both $\xi$ and $-K_{\Pee(V)}$ are ample divisors.
Thus, $\Pee(V)$ is a Fano variety, and its quantum cohomology ring
is well-defined \cite{13}.  Here we choose the symplectic form $\omega$ on
$\Pee(V)$ to be the Kahler form $\omega$ such that $[\omega] = -K_{\Pee(V)}$.
Let $f_\omega^1$ and $f_\omega^2$ be the two relations in (1.1)
evaluated in the quantum cohomology ring
$H^*_{\omega}(\Pee(V); \Zee)$. Then by the Theorem 2.2 in \cite{14},
the quantum cohomology $H^*_{\omega}(\Pee(V); \Zee)$ is
the ring generated by $h$ and $\xi$ with the two relations
$f_\omega^1$ and $f_\omega^2$:
$$H_\omega^*(\Pee(V); \Zee) = \Zee [h, \xi]/(f_\omega^1, f_\omega^2)
\eqno (1.3)$$

By Mori's Cone Theorem \cite{5}, $\Pee(V)$ has exactly
two extremal rays $R_1$ and $R_2$. Up to an order of $R_1$ and $R_2$,
the integral generator $A_1$ of $R_1$ is represented by lines
in the fibers of the projection $\pi$.
We shall show that under certain numerical conditions,
the nontrivial homology classes
$A \in H_2(\Pee(V); \Zee)$ which give nontrivial quantum corrections
are $A_1$ and $A_2$, where $A_2$ is represented by
some smooth rational curves in $\Pee(V)$ which are isomorphic
to lines in $\Pee^n$ via $\pi$. In general, it is unclear whether
$A_2$ generates the second extremal ray $R_2$. However, we shall prove
that under further restrictions on $V$, $A_2$ generates
the extremal ray $R_2$. These analyses enable us to determine
the quantum cohomology ring $H^*_{\omega}(\Pee(V); \Zee)$.

The simplest ample bundle over $\Pee^n$ is perhaps the direct sum of
line bundles $V=\oplus^r_{i=1} {\Cal O}(m_i)$ where $m_i>0$ for every $i$.
Since we can twist $V$ by ${\Cal O}(-1)$ without changing $\Pee(V)$,
we can assume that $\hbox{min}\{ m_1, \ldots, m_r \}=1$.
In this case, $\Pee(V)$ is a special case of toric variety.
Batyrev \cite{2} conjectured a general formula for quantum cohomology of
toric varieties. Furthermore, he computed the contributions
from certain moduli spaces of holomorphic maps
which have expected dimensions. In our case,
the contributions Batyrev computed are only part of the data
to compute the quantum cohomology. As we explained earlier,
the difficulty in our case lies precisely in computing the contributions
from the moduli spaces with wrong dimensions. Nevertheless, in our case,
Batyrev's formula (see also \cite{1}) reads as follows.
\vskip 0.1in
\noindent
{\bf Batyrev's Conjecture:}
{\it Let $V=\oplus^r_{i=1} {\Cal O}(m_i)$ where $m_i > 0$ for every $i$.
Then the quantum cohomology ring $H^*_{\omega}(\Pee(V); \Zee)$ is
generated by $h$ and $\xi$ with two relations}
$$h^{n+1}=\prod^r_{i=1}(\xi-m_ih)^{m_i-1} \cdot e^{-t(n+1+r-\sum_{i=1}^r m_i)}
\qquad {and} \qquad
\prod^r_{i=1}(\xi-m_ih) = e^{-tr}.$$
\vskip 0.1in

Our first result partially verifies Batyrev's conjecture.

\theorem{A} Batyrev's conjecture holds if
$$\sum_{i=1}^r m_i < \text{min}(2r, (n+1+2r)/2, (2n+2+r)/2).$$
\endproclaim

Note that under the numerical condition of Theorem A,
only extremal rational curves with fundamental classes $A_1$ and $A_2$
give the contributions to the two relations in the quantum cohomology.
The moduli space of rational curves $\frak M(A_2, 0)$ with fundamental
class $A_2$ does not have expected dimension in general. But it is compact.
This fact simplifies a great deal of the excessive intersection theory
involved.
To remove the numerical condition, we have to consider other moduli spaces
(for example $\frak M(kA_2, 0)$ with $k > 1$ and
its excessive intersection theory).
These moduli spaces are not compact in general.
Then, we have an extra difficulty of the compactification and
the appropriate excessive intersection theory with it.
It seems to be a difficult problem and we shall not pursue here.

In general, ample bundles over $\Pee^n$ are not direct sums of line bundles.
We can say much less about its quantum cohomology. However,
we obtain some result about its general form and
compute the leading coefficient.

\theorem{B} {\rm (i)} Let $V$ be a rank-$r$ ample bundle over $\Pee^n$.
Assume either $c_1 \le n$ or $c_1 \le (n + r)$ and
$V\otimes {\Cal O}_{\Pee^n}(-1)$ is nef so that $\Pee(V)$ is Fano.
Then the quantum cohomology $H^*_{\omega}(\Pee(V); \Zee)$ is
the ring generated by $h$ and $\xi$ with two relations
$$h^{n+1} = \sum_{i+j \le (c_1 - r)} a_{i,j} \cdot h^i \cdot \xi^j
\cdot e^{-t(n+1-i-j)}$$
$$\sum_{i=0}^r (-1)^i c_i \cdot h^i \cdot \xi^{r-i} = e^{-tr}
+ \sum_{i+j \le (c_1-n-1)} b_{i,j} \cdot h^i \cdot \xi^{j}
\cdot e^{-t(r-i-j)}$$
where the coefficients $a_{i, j}$ and $b_{i, j}$ are integers depending on $V$;

{\rm (ii)} If we further assume that $c_1<2r$,
then the leading coefficient $a_{0,c_1-r}=1$.
\endproclaim

It is understood that when $c_1 \le n$, then the summation
$\sum_{i+j \le (c_1-n-1)}$
in the second relation in Theorem B (i) does not exist.
In general, it is not easy to determine all the integers $a_{i, j}$ and
$b_{i, j}$ in Theorem B (i). However, it is possible to
compute these numbers when $(c_1 - r)$ is relatively small. For instance,
when $(c_1 - r) = 0$, then necessarily $V = \Cal O_{\Pee^n}(1)^{\oplus r}$
and it is well-known that the quantum cohomology
$H^*_{\omega}(\Pee(V); \Zee)$ is the ring generated by $h$ and $\xi$
with the two relations $h^{n+1} = e^{-t(n+1)}$ and
$\sum_{i=0}^r (-1)^i c_i \cdot h^i \cdot \xi^{r-i} = e^{-tr}$.
When $(c_1 - r) = 1$ and $r < n$,
then necessarily $V = \Cal O_{\Pee^n}(1)^{\oplus (r - 1)}
\oplus \Cal O_{\Pee^n}(2)$.  When $(c_1 - r) = 1$ and $r = n$,
then $V = \Cal O_{\Pee^n}(1)^{\oplus (r - 1)} \oplus \Cal O_{\Pee^n}(2)$
or $V = T_{\Pee^n}$ the tangent bundle of $\Pee^n$.
In these cases, $V\otimes {\Cal O}_{\Pee^n}(-1)$ is nef. In particular, the
direct sum cases have been computed by Theorem A.
We shall  prove the following.

\proposition{C} The quantum cohomology ring
$H^*_{\omega}(\Pee(T_{\Pee^n}); \Zee)$ with $n \ge 2$ is the ring
generated by $h$ and $\xi$ with the two relations:
$$h^{n+1} = \xi \cdot e^{-tn} \qquad \text{and}
\qquad \sum_{i=0}^n (-1)^i c_i \cdot h^i \cdot \xi^{n-i}
= (1 + (-1)^n) \cdot e^{-tn}.$$
\endproclaim

Recall that for an arbitrary projective bundle over a general manifold,
its cohomology ring is a module over the cohomology ring of the base
with the generator $\xi$ and the second relation of (1.1).
Naively. one may think that the quantum cohomology of projective bundle
is a module over the quantum cohomology of base with the generator $\xi$
and the quantanized second relation. Our calculation shows that
one can not expect such simplicity for its quantum cohomology ring.
We hope that our results could shed some light
on the quantum cohomology for general projective bundles,
which we shall leave for future research.

Our paper is organized as follows.
In section 2, we discuss the extremal rays and extremal rational curves.
In section 3, we review the definition of quantum product
and compute some Gromov-Witten invariants.
In the remaining three sections, we prove Theorem B, Theorem A,
and Proposition C respectively.

\medskip\noindent
{\bf Acknowledgements:} We would like to thank Sheldon Katz, Yungang Ye,
and Qi Zhang for valuable helps and stimulating discussions.
In particular, we are grateful to Sheldon Katz for bringing us
the attention of Batyrev's conjecture.

\section{2. Extremal rational curves}

Assume that $V$ is ample such that either $c_1 \le (n+1)$
or $c_1 \le (n + r)$ and $V\otimes {\Cal O}_{\Pee^n}(-1)$ is nef.
In this section, we study the extremal rays and extremal rational curves
in the Fano variety $\Pee(V)$. By Mori's Cone Theorem (p.25 in \cite{5}),
$\Pee(V)$ has precisely two extremal rays $R_1 = \Bbb R_{\ge 0} \cdot A_1$
and $R_2 = \Bbb R_{\ge 0} \cdot A_2$ such that the cone
$\hbox{NE}(\Pee(V))$ of curves in $\Pee(V)$ is equal to $R_1 + R_2$ and
that $A_1$ and $A_2$ are the homology classes of two rational curves
$E_1$ and $E_2$ in $\Pee(V)$ with
$0 < -K_{\Pee(V)}(A_i) \le \hbox{dim}(\Pee(V)) + 1$.
Up to orders of $A_1$ and $A_2$, we have $A_1 = (h_{n}\xi_{r-2})_*$,
that is, $A_1$ is represented by lines in the fibers of $\pi$.
It is also well-known that if $V = \oplus_{i=1}^r \Cal O_{\Pee^n}(m_i)$
with $m_1 \le \ldots \le m_r$,
then $A_2 = [h_{n-1}\xi_{r-1}+ (m_1 - c_1) h_{n}\xi_{r-2}]_*$
which is represented by a smooth rational curve in $\Pee(V)$
isomorphic to a line in $\Pee^n$ via $\pi$. However, in general,
it is not easy to determine the homology class $A_2$ and
the extremal rational curves representing $A_2$. Assume that
$$V|_\ell = \oplus_{i = 1}^r \Cal O_\ell(m_i) \eqno (2.1)$$
for generic lines $\ell \subset \Pee^n$ where we let $m_1 \le \ldots \le m_r$.
Since $V$ is ample, $m_1 \ge 1$.

\lemma{2.2} Let $A = [h_{n-1}\xi_{r-1} + (m_1 - c_1) h_{n}\xi_{r-2}]_*$. Then,

\roster
\item"{(i)}" $A$ is represented by a smooth rational curve isomorphic to
a line in $\Pee^n$;

\item"{(ii)}" $A_2 = A$ if and only if $(\xi- m_1 h)$ is nef;

\item"{(iii)}" $A_2 = A$ if $2c_1 \le (n + 1)$;

\item"{(iv)}" $A$ can not be represented by reducible or nonreduced curves
if $m_1 = 1$.
\endroster
\endproclaim
\proof (i) Let $\ell \subset \Pee^n$ be a generic line.
Then we have a natural projection $V|_\ell = \oplus_{i = 1}^r \Cal O_\ell(m_i)
\to \Cal O_\ell(m_1)$. By the Proposition 7.12 in Chapter II of \cite{7},
this surjective map $V|_\ell \to \Cal O_\ell(m_1) \to 0$ induces
a morphism $g: \ell \to \Pee(V)$. Then $g(\ell)$ is isomorphic to $\ell$
via the projection $\pi$. Since $h([g(\ell)]) = 1$ and $\xi([g(\ell)]) = m_1$,
we have
$$[g(\ell)] = [h_{n-1}\xi_{r-1}+ (m_1 - c_1) h_{n}\xi_{r-2}]_* = A.$$

(ii) First of all, if
$A_2 = [h_{n-1}\xi_{r-1}+ (m_1 - c_1) h_{n}\xi_{r-2}]_*$,
then for any curve $E$, $[E] = a(h_{n}\xi_{r-2})_* +
b[h_{n-1}\xi_{r-1}+ (m_1 - c_1) h_{n}\xi_{r-2}]_*$ for some nonnegative numbers
$a$ and $b$; so $(\xi- m_1h)([E]) = a \ge 0$; therefore $(\xi- m_1h)$ is nef.
Conversely, if $(\xi- m_1h)$ is nef, then $0 \le (\xi- m_1h)([E])
= ac_1 + b - am_1$ where $[E] = (ah_{n-1}\xi_{r-1}+ b h_{n}\xi_{r-2})_*$
for some curve $E$; thus $[E] = (ac_1 + b - am_1)(h_{n}\xi_{r-2})_* +
a[h_{n-1}\xi_{r-1}+ (m_1 - c_1) h_{n}\xi_{r-2}]_*$;
it follows that $A_2 = [h_{n-1}\xi_{r-1}+ (m_1 - c_1) h_{n}\xi_{r-2}]_* = A$.

(iii) Let $A_2 = (ah_{n-1}\xi_{r-1} + bh_{n}\xi_{r-2})_*$.
Since $A_1 = (h_{n}\xi_{r-2})_*$ and $a = h(A_2) \ge 0$, $a \ge 1$.
If $a > 1$, then since $2c_1 \le (n + 1)$, we see that
$$\align
-K_{\Pee(V)}(A_2) &= (n+1-c_1)a +r \cdot \xi(A_2) \ge 2(n+1-c_1) + r \\
&> n + r = \text{dim}(\Pee(V)) + 1; \\
\endalign$$
but this contradicts with $-K_{\Pee(V)}(A_2) \le \hbox{dim}(\Pee(V)) + 1$.
Thus $a = 1$ and $A_2 = (h_{n-1}\xi_{r-1} + bh_{n}\xi_{r-2})_*$.
Now $[\pi(E_2)] = \pi_*(A_2) = (h_{n-1})_*$.
So $\pi(E_2)$ is a line in $\Pee^n$.
Since $V|_\ell = \oplus_{i = 1}^r \Cal O_\ell(m_i)$ for a generic line
$\ell \subset \Pee^n$,
$V|_{\pi(E_2)} = \oplus_{i = 1}^r \Cal O_{\pi(E_2)}(m_i')$
where $m_i' \ge m_1$ for every $i$. Thus, $\xi(A_2) \ge m_1$,
and so $c_1 + b \ge m_1$. It follows that
$$A_2 = [h_{n-1}\xi_{r-1} + (m_1- c_1) h_{n}\xi_{r-2}]_*
+ (c_1 + b - m_1) \cdot (h_{n}\xi_{r-2})_*.$$
Therefore, $A_2 = [h_{n-1}\xi_{r-1}+ (m_1 - c_1) h_{n}\xi_{r-2}]_* = A$.

(iv) Since $\xi(A) = m_1 = 1$ and $\xi$ is ample, the conclusion follows.
\endproof

Next, let $\frak M(A,0)$ be the moduli space of morphisms
$f: \Pee^1 \to \Pee(V)$ with $[\hbox{Im}(f)] = A$.
In the lemma below, we study the morphisms in $\frak M(A,0)$
when $A = [h_{n-1}\xi_{r-1}+ (m - c_1) h_{n}\xi_{r-2}]_*$.
Note that $\xi(A) = m$.

\lemma{2.3} Let $A = [h_{n-1}\xi_{r-1}+ (m - c_1) h_{n}\xi_{r-2}]_*$.

\roster
\item"{(i)}" If $\frak M(A, 0) \ne \emptyset$, then $m \ge m_1$ and
$\frak M(A, 0)$ consists of embeddings $f: \ell \to \Pee(V)$
induced by surjective maps $V|_\ell \to \Cal O_\ell(m) \to 0$
where $\ell$ are lines in $\Pee^n$;

\item"{(ii)}" If $m = m_1$ and $m_1 = \ldots = m_k < m_{k+1} \le
\ldots \le m_r$, then the moduli space $\frak M(A, 0)$
has (complex) dimension $(2n +k)$;

\item"{(iii)}" If $m \ge m_r$, then $\frak M(A, 0)$ has dimension
$(2n+r+rm-c_1)$.
\endroster
\endproclaim
\proof
(i) Let $f: \Pee^1 \to \Pee(V)$ be a morphism in $\frak M(A, 0)$.
Then $[\hbox{Im}(f)] = A = [h_{n-1}\xi_{r-1}+ (m - c_1) h_{n}\xi_{r-2}]_*$.
Since $h(A) = 1$, $\pi^*H \cap f(\Pee^1)$ consists of a single point
for any hyperplane $H$ in $\Pee^n$. Thus,
$\pi|_{f(\Pee^1)}: f(\Pee^1) \to (\pi \circ f)(\Pee^1)$ is an isomorphism
and $\ell = (\pi \circ f)(\Pee^1)$ is a line in $\Pee^n$.
Since $h([\ell]) = 1$, $(\pi \circ f): \Pee^1 \to \ell = (\pi \circ f)(\Pee^1)$
is also an isomorphism, and so is $f: \Pee^1 \to f(\Pee^1)$.
Replacing $f: \Pee^1 \to \Pee(V)$ by
$f \circ (\pi \circ f)^{-1}: \ell \to \Pee(V)$,
we conclude that $\frak M(A, 0)$ consists of embeddings $f: \ell \to \Pee(V)$
such that $[\hbox{Im}(f)] = A$, $\ell$ are lines in $\Pee^n$,
and $\pi|_{f(\ell)}: f(\ell) \to \ell$ are isomorphisms.
In particular, these embeddings $f: \ell \to \Pee(V)$ are
sections to the natural projection
$\pi|_{\Pee(V|_\ell)}: \Pee(V|_\ell) \to \ell$.
Thus, by the Proposition 7.12 in Chapter II of \cite{7},
these embeddings are induced by surjective maps
$V|_\ell \to \Cal O_\ell(m) \to 0$.
By (2.1), the splitting type of the restrictions of $V$ to generic lines
in $\Pee^n$ is $(m_1, \ldots, m_r)$ with $m_1 \le \ldots \le m_r$;
thus we must have $V|_\ell = \oplus_{i = 1}^r \Cal O_\ell(m_i')$ where
$m_i' \ge m_1$ for every $i$. It follows that
$m \ge \hbox{min} \{ m_1', \ldots, m_r' \} \ge m_1$.

(ii) Note that all the lines in $\Pee^n$ are parameterized by
the Grassmannian $G(2, n+1)$ which has dimension $2(n -1)$.
For a fixed generic line $\ell \subset \Pee^n$, the surjective maps
$V|_\ell \to \Cal O_\ell(m_1) \to 0$ are parameterized by
$$\Pee(\hbox{Hom}(V|_\ell, \Cal O_\ell(m_1))) \cong
\Pee(\oplus_{i=1}^r H^0(\ell, \Cal O_\ell(m_1-m_i))) \cong \Pee^{k - 1};$$
It follows from (i) that as the generic line $\ell$ varies,
the morphisms $f: \ell \to \Pee(V)$ induced by these surjective maps
$V|_\ell \to \Cal O_\ell(m_1) \to 0$ form an open dense subset of
$\frak M(A, 0)$. Thus, $\hbox{dim}(\frak M(A, 0)) = 3+ 2(n -1) + (k - 1)
= 2n +k$.

(iii) As in the proof of (ii), for a fixed generic line $\ell \subset \Pee^n$,
the surjective maps $V|_\ell \to \Cal O_\ell(m) \to 0$ are
parameterized by a nonempty open subset of
$$\Pee(\hbox{Hom}(V|_\ell, \Cal O_\ell(m))) \cong
\Pee(\oplus_{i=1}^r H^0(\ell, \Cal O_\ell(m-m_i)))
\cong \Pee^{(rm - c_1 +r) - 1}.$$
As the generic line $\ell$ varies, the morphisms $f: \ell \to \Pee(V)$
induced by these surjective maps $V|_\ell \to \Cal O_\ell(m) \to 0$
form an open dense subset of $\frak M(A, 0)$. It follows that
$\frak M(A, 0)$ has dimension $(2n+r+rm-c_1)$.
\endproof

\section{3. Calculation of Gromov-Witten invariants}

In this section, we shall compute some Gromov-Witten invariants of $\Pee(V)$.
First of all, we recall that for two homogeneous elements
$\alpha$ and $\beta$ in $H^*(\Pee(V); \Zee)$, the quantum product
$\alpha \cdot \beta \in H^*(\Pee(V); \Zee)$ can be written as
$$\alpha \cdot \beta = \sum_{A \in H_2(\Pee(V); \Zee)}
(\alpha \cdot \beta)_A \cdot e^{t \cdot K_{\Pee(V)}(A)} \eqno (3.1)$$
where $(\alpha \cdot \beta)_A$ has degree
$\hbox{deg}(\alpha) + \hbox{deg}(\beta) + 2K_{\Pee(V)}(A)$
and is defined by
$$(\alpha \cdot \beta)_A(\gamma_*) = \Phi_{(A, 0)}(\alpha, \beta, \gamma)$$
for a homogeneous cohomology class $\gamma \in H^*(\Pee(V); \Zee)$ with
$$\hbox{deg}(\gamma) = -2K_{\Pee(V)}(A) + 2(n+r-1) - \hbox{deg}(\alpha)
- \hbox{deg}(\beta). \eqno (3.2)$$
Furthermore, for higher quantum products, we have
$$\alpha_1 \cdot \alpha_2 \cdot \ldots \cdot \alpha_k =
\sum_{A \in H_2(\Pee(V); \Zee)} (\alpha_1 \cdot \alpha_2 \cdot
\ldots \cdot \alpha_k)_A \cdot e^{t \cdot K_{\Pee(V)}(A)} \eqno (3.3)$$
where $(\alpha_1 \cdot \alpha_2 \cdot \ldots \cdot \alpha_k)_A$ is defined as
$(\alpha_1 \cdot \alpha_2 \cdot \ldots \cdot \alpha_k)_A(\gamma_*)
= \Phi_{(A,0)}(\alpha_1, \alpha_2, \ldots, \alpha_k, \gamma)$. Thus,
$\alpha_1 \cdot \alpha_2 \cdot \cdots \cdot \alpha_k
= \alpha_1\alpha_2\dots\alpha_k + \text{(lower\ order\ terms)}$,
where $\alpha_1\alpha_2\dots\alpha_k$
stands for the ordinary cohomology product of
$\alpha_1, \alpha_2, \ldots, \alpha_k$,
and the degree of a lower order term is dropped by $2K_{\Pee(V)}(A)$
for some $A \in H_2(\Pee(V); \Zee)$ which is represented by
a nonconstant effective rational curve.

There are two explanations for the Gromov-Witten invariant
$\Phi_{(A, 0)}(\alpha, \beta, \gamma)$ defined by the second author \cite{12}.
Recall that the Gromov-Witten invariant is only defined for
a generic almost complex structure and that $\frak M(A,0)$ is
the moduli space of morphisms $f: \Pee^1 \to \Pee(V)$
with $[\hbox{Im}(f)] = A$. Assume the genericity conditions:
\roster
\item"{(i)}" $\frak M(A, 0)/PSL(2; \Cee)$ is smooth
in the sense that $h^1(N_f) = 0$ for every $f \in \frak M(A, 0)$ where
$N_f$ is the normal bundle, and
\item"{(ii)}"  the homology class $A$ is only
represented by irreducible and reduced curves.
\endroster
\noindent
Then the complex structure is already generic and one can use
algebraic geometry to calculate the Gromov-Witten invariants.
Moreover, $\frak M(A, 0)/PSL(2; \Cee)$ is
compact with the expected complex dimension
$$-K_{\Pee(V)}(A) + (n+r-1) - 3. \eqno (3.4)$$
The first explanation for $\Phi_{(A, 0)}(\alpha, \beta, \gamma)$ is that
when $\alpha, \beta, \gamma$ are classes of subvarieties $B, C, D$
of $\Pee(V)$ in general position, $\Phi_{(A, 0)}(\alpha, \beta, \gamma)$
is the number of rational curves $E$ in $\Pee(V)$ such that
$[E] = A$ and that $E$ intersects with $B, C, D$
(counted with suitable multiplicity). The second explanation for
$\Phi_{(A, 0)}(\alpha, \beta, \gamma)$ is that
$$\Phi_{(A, 0)}(\alpha, \beta, \gamma) = \int_{\frak M(A, 0)} e_0^*(\alpha)
\cdot e_1^*(\beta) \cdot e_2^*(\gamma) $$
where the evaluation map $e_i: \frak M(A, 0) \to \Pee(V)$ is defined by
$e_i(f) = f(i)$.

Assume that the genericity condition (i) is not satisfied but
$h^1(N_f)$ is independent of $f \in \frak M(A, 0)$ and
$\frak M(A, 0)/PSL(2; \Cee)$ is smooth with dimension
$$-K_{\Pee(V)}(A) + (n+r-1) - 3+h^1(N_f).$$
Then one can form an obstruction bundle $COB$ of rank $h^1(N_f)$
over the moduli space $\frak M(A, 0)$. Moreover,
if the genericity condition (ii) is satisfied,
then by the Proposition 5.7 in \cite{11}, we have
$$\Phi_{(A, 0)}(\alpha, \beta, \gamma) = \int_{\frak M(A, 0)} e_0^*(\alpha)
\cdot e_1^*(\beta) \cdot e_2^*(\gamma) \cdot e(COB) \eqno (3.5)$$
where $e(COB)$ stands for the Euler class of the bundle $COB$.

We remark that in general, the cohomology class $h_i \xi_j$ may not
be able to be represented by a subvariety of $\Pee(V)$.
However, since $\xi$ is ample, $s \xi$ is very ample for $s \gg 0$.
Thus, the multiple $t h_i \xi_j$ with $t \gg 0$ can be represented
by a subvariety of $\Pee(V)$ whose image in $\Pee^n$ is
a linear subspace of codimension $i$.
Since $\Phi_{(A, 0)}(\alpha, \beta, h_i \xi_j) =
{1/t} \cdot \Phi_{(A, 0)}(\alpha, \beta, t \cdot h_i \xi_j)$
for $\alpha$ and $\beta$ in $H^*(\Pee(V); \Zee)$,
it follows that to compute $\Phi_{(A, 0)}(\alpha, \beta, h_i \xi_j)$,
it suffices to compute $\Phi_{(A, 0)}(\alpha, \beta, t \cdot h_i \xi_j)$.
In the proofs below, we shall assume implicitly that $t = 1$ for simplicity.

Now we compute the Gromov-Witten invariant
$\Phi_{((h_{n}\xi_{r-2})_*, 0)}(\xi, \xi_{r-1}, h_{n}\xi_{r-1})$.

\lemma{3.6}
$\Phi_{((h_{n}\xi_{r-2})_*, 0)}(\xi, \xi_{r-1}, h_{n}\xi_{r-1}) = 1$.
\endproclaim
\proof
First of all, we notice that $A = (h_{n} \xi_{r-2})_*$ can only be
represented by lines $\ell$ in the fibers of $\pi$. In particular,
there is no reducible or nonreduced effective curves representing $A$.
Thus, $\frak M(A, 0)/PSL(2; \Cee)$ is compact and has dimension:
$$\text{dim}(\Pee^n) + \text{dim} G(2, r)= n + 2(r-2) = n + 2r - 4$$
which is the expected dimension by (3.4) (here we use $G(2, r)$ to stand for
the Grassmannian of lines in $\Pee^{r-1}$). Next, we want to show that
$\frak M(A, 0)/PSL(2; \Cee)$ is smooth. Let $p = \pi(\ell)$.
Then from the two inclusions $\ell \subset \pi^{-1}(p) \subset \Pee(V)$,
we obtain an exact sequence relating normal bundles:
$$0 \to N_{\ell|\pi^{-1}(p)} \to N_{\ell|\Pee(V)} \to
(N_{\pi^{-1}(p)|\Pee(V)})|_\ell \to 0.$$
Since $N_{\ell|\pi^{-1}(p)} = N_{\ell|\Pee^{r-1}}
= \Cal O_{\ell}(1)^{\oplus (r-2)}$
and $N_{\pi^{-1}(p)|\Pee(V)} = (\pi|_{\pi^{-1}(p)})^*T_{p, \Pee^n}$,
the previous exact sequence is simplified into the exact sequence
$$0 \to \Cal O_{\ell}(1)^{\oplus (r-2)} \to N_{\ell|\Pee(V)} \to
(\pi|_\ell)^*T_{p, \Pee^n} \to 0.$$
It follows that $H^1(\ell, N_{\ell|\Pee(V)}) = 0$. Thus,
$\frak M(A, 0)/PSL(2; \Cee)$ is smooth.

Finally, the Poincar\'e dual of $h_n \xi_{r-1}$ is represented by a point
$q_0 \in \Pee(V)$. If a line $\ell \in \frak M(A,0)$ intersects $q_0$,
then $\ell \subset \pi^{-1}(\pi(q_0))$. Since the restriction of $\xi$ to
the fiber $\pi^{-1}(\pi(q_0)) \cong \Pee^{r-1}$ is
the cohomology class of a hyperplane in $\Pee^{r-1}$, we conclude that
$\Phi_{((h_{n}\xi_{r-2})_*, 0)}(\xi, \xi_{r-1}, h_{n}\xi_{r-1}) = 1$.
\endproof

Next we show the vanishing of some Gromov-Witten invariant.

\lemma{3.7} Let $A = b(h_n \xi_{r-2})_*$ with $b \ge 1$
and $\alpha \in H^*(\Pee(V); \Zee)$. Then,
$$\Phi_{(A, 0)}(h_{p_1}\xi_{q_1}, h_{p_2}\xi_{q_2}, \alpha) = 0$$
if $p_1, q_1, p_2, q_2$ are nonnegative integers with $(q_1+q_2) < r$.
\endproclaim
\proof
We may assume that $\alpha$ is a homogeneous class in
$H^*(\Pee(V); \Zee)$. By (3.2),
$$\align
{1 \over 2} \cdot \hbox{deg}(\alpha)
&= (n+r-1) - K_{\Pee(V)}(A) - (p_1+p_2+q_1+q_2) \\
&= (n+r+br-1)-(p_1+p_2+q_1+q_2). \\
\endalign$$
Let $\alpha = h_{(n+r+br-1)-(p_1+p_2+q_1+q_2+q_3)} \xi_{q_3}$ with
$0 \le q_3 \le (r - 1)$.
Let $B, C, D$ be the subvarieties of $\Pee(V)$ in general position,
whose homology classes are Poincar\'e dual to
$h_{p_1}\xi_{q_1}, h_{p_2}\xi_{q_2}, \alpha$ respectively.
Then the homology classes of $\pi(B), \pi(C), \pi(D)$ in $\Pee^n$ are
Poincar\'e dual to $h_{p_1}, h_{p_2},
h_{(n+r+br-1)-(p_1+p_2+q_1+q_2+q_3)}$ respectively.
Since $(q_1+q_2+q_3) < (2r -1)$, we have
$p_1 + p_2 + [(n+r+br-1)-(p_1+p_2+q_1+q_2+q_3)]
= (n+r+br-1)-(q_1+q_2+q_3) > n$.
Thus, $\pi(B) \cap \pi(C) \cap \pi(D) = \emptyset$.
Notice that the genericity conditions (i) and (ii)
mentioned earlier in this section are not satisfied for $b\ge 2$.
However, we observe that these conditions can be relaxed by assuming:
\roster
\item"{(i$'$)}" $h^1(N_f)=0$ for every $f\in \frak M(A, 0)$ such that
$\hbox{Im}(f)$ intersects $B, C, D$, and
\item"{(ii$'$)}" there is no reducible or nonreduced effective
(connected) curve $E$ such that $[E]=A$ and $E$ intersects $B, C, D$.
\endroster
\noindent
In fact, we will show that there is no effective connected curve $E$ at all
representing $A$ and intersecting $B, C, D$.
It obviously implies (i$'$), (ii$'$) and
$$\Phi_{(A, 0)}(h_{p_1}\xi_{q_1}, h_{p_2}\xi_{q_2}, \alpha) = 0.$$
Suppose that $E=\sum a_i E_i$ is such an effective connected curve
where $a_i > 0$ and $E_i$ is irreducible and reduced.
Then, $\sum a_i [E_i] = [E] = A$.
Since $(h_n \xi_{r-2})_*$ generates an extremal ray for $\Pee(V)$,
$[E_i]=b_i (h_n \xi_{r-2})_*$ for $0<b_i\le b$.
Thus the curves $E_i$ are contained in the fibers of $\pi$.
Since $E$ is connected, all the curves $E_i$ must be contained
in the same fiber of $\pi$. So $\pi(E)$ is a single point.
Since $E$ intersects $B, C, D$, $\pi(E)$ intersects with
$\pi(B), \pi(C), \pi(D)$.
It follows that $\pi(B) \cap \pi(C) \cap \pi(D)$ contains $\pi(E)$
and is nonempty. Therefore we obtain a contradiction.
\endproof

Finally, we show that if $c_1 < 2r$ and
$A = [h_{n-1}\xi_{r-1} + (1 - c_1) h_{n}\xi_{r-2}]_*$,
then $\Phi_{(A, 0)}(h, h_n, h_n \xi_{2r-c_1-1}) = 1$.
Since $c_1 < 2r$, we see that for a generic line $\ell \subset \Pee^n$,
$$V|_{\ell} = \Cal O_{\ell}(1)^{\oplus k}
\oplus \Cal O_{\ell}(m_{k+1}) \oplus \ldots \oplus \Cal O_{\ell}(m_r)$$
where $k \ge 1$ and $2 \le m_{k+1} \le \ldots \le m_r$.
We remark that even though the moduli space $\frak M(A, 0)/PSL(2; \Cee)$
is compact by Lemma 2.2 (iv),
it may not have the correct dimension by Lemma 2.3 (ii).
The proof is lengthy, but the basic idea is that
we shall determine the obstruction bundle and
use the formula (3.5).

\lemma{3.8} Let $V$ be a rank-$r$ ample vector bundle over
$\Pee^n$ satisfying $c_1<2r$ and the assumption of Theorem B (i).
If $A = [h_{n-1}\xi_{r-1} + (1 - c_1) h_{n}\xi_{r-2}]_*$, then
$$\Phi_{(A, 0)}(h, h_n, h_n \xi_{2r-c_1-1}) = 1.$$
\endproclaim
\noindent
{\it Proof.} Note that by Lemma 2.2 (iv),
the moduli space $\frak M(A, 0)/PSL(2; \Cee)$ is compact.
Let $B, C, D$ be the subvarieties of $\Pee(V)$ in general position,
whose homology classes are Poincar\'e dual to
$h, h_n, h_n \xi_{2r-c_1-1}$ respectively. Then the homology classes of
$\pi(B), \pi(C)$, $\pi(D)$ in $\Pee^n$ are Poincar\'e dual to
$h, h_n, h_n$ respectively. Thus $\pi(C)$ and $\pi(D)$ are
two different points in $\Pee^n$. Let $\ell_0$ be the unique line passing
$\pi(C)$ and $\pi(D)$. Let $V|_{\ell_0} = \Cal O_{\ell_0}(1)^{\oplus k}
\oplus \Cal O_{\ell_0}(m_{k+1}) \oplus \ldots \oplus \Cal O_{\ell_0}(m_r)$
where $2 \le m_{k+1} \le \ldots \le m_r$. Since $c_1 < 2r$, $k \ge 1$.
Let $f: \ell \to \Pee(V)$ be
a morphism in $\frak M(A, 0)$ for some line $\ell \in \Pee^n$.
If $\hbox{Im}(f)$ intersects with $B, C$, and $D$, then $\ell = \ell_0$.
As in the proof of Lemma 2.3 (ii), the morphisms $f: \ell_0 \to \Pee(V)$
in $\frak M(A, 0)$ are parameterized by
$\Pee(\hbox{Hom}(V|_{\ell_0}, \Cal O_{\ell_0}(1))) \cong \Pee^{k - 1}$;
moreover, $\hbox{Im}(f)$ are of the form:
$$\ell_0 \times \{ q \} \subset \ell_0 \times \Pee^{k - 1}
= \Pee(\Cal O_{\ell_0}(1)^{\oplus k})
\subset \Pee(V|_{\ell_0}) \subset \Pee(V) \eqno (3.9)$$
where $q$ stands for points in $\Pee^{k - 1} \subset \Pee^{r - 1} \cong
\pi^{-1}(\pi(D))$. Note that $\ell_0 \times \{ q \}$ always intersects
with $B$ and $C$, and that $D$ is a dimension-$(c_1 - r)$ linear subspace
in $\Pee^{r - 1} \cong \pi^{-1}(\pi(D))$. Thus, $\ell_0 \times \{ q \}$
intersects with $B, C, D$ simultaneously if and only if
$\ell_0 \times \{ q \}$ intersects with $D$, and if only only if
$$q \in \Pee^{c_1+k-2r} \overset \hbox{def} \to =
\Pee^{k - 1} \cap D \subset \Pee^{r - 1} \cong \pi^{-1}(\pi(D)).
\eqno (3.10)$$
It follows that $\frak M/PSL(2; \Cee) \cong \Pee^{c_1+k-2r}$ where
$\frak M$ consists of morphisms $f \in \frak M(A, 0)$ such that
$\hbox{Im}(f)$ intersects with $B, C, D$ simultaneously.

If $c_1+k-2r = 0$, then $a_0 = \Phi_{(A, 0)}(h, h_n, h_n \xi_{2r-c_1-1}) = 1$.
But in general, we have $c_1+k-2r \ge 0$. We shall use (3.5) to compute
$a_0 = \Phi_{(A, 0)}(h, h_n, h_n \xi_{2r-c_1-1})$.
Let $N_f = N_{\ell_0 \times \{ q \}|\Pee(V)}$ be the normal bundle
of $\hbox{Im}(f) = \ell_0 \times \{ q \}$ in $\Pee(V)$.
If $h^1(N_f)$ is constant for every $f \in \frak M$,
then by (3.5), $\Phi_{(A, 0)}(h, h_n, h_n \xi_{2r-c_1-1})$
is the Euler number $e(COB)$ of the rank-$(c_1+k-2r)$ obstruction bundle
$COB$ over
$$\frak M/PSL(2; \Cee) \cong \Pee^{c_1+k-2r}.$$
Thus we need to show that $h^1(N_f)$ is constant for every $f \in \frak M$.

First, we study the normal bundle
$N_{\ell_0 \times \Pee^{c_1+k-2r}|\Pee(V)}$. The three inclusions
$$\ell_0 \times \Pee^{c_1+k-2r}
\subset \ell_0 \times \Pee^{k - 1} = \Pee(\Cal O_{\ell_0}(1)^{\oplus k})
\subset \Pee(V|_{\ell_0}) \subset \Pee(V) \eqno (3.11)$$
give rise to two exact sequences relating normal bundles:
$$0 \to N_{\ell_0 \times \Pee^{c_1+k-2r}|\Pee(V|_{\ell_0})}
\to N_{\ell_0 \times \Pee^{c_1+k-2r}|\Pee(V)} \to
N_{\Pee(V|_{\ell_0})|\Pee(V)} \to 0$$
$$0 \to N_{\ell_0 \times \Pee^{c_1+k-2r}|\Pee(\Cal O_{\ell_0}(1)^{\oplus k})}
\to N_{\ell_0 \times \Pee^{c_1+k-2r}|\Pee(V|_{\ell_0})} \to
N_{\Pee(\Cal O_{\ell_0}(1)^{\oplus k})|\Pee(V|_{\ell_0})} \to 0$$
Notice that $N_{\Pee(V|_{\ell_0}) | \Pee(V)} =
(\pi|_{\Pee(V|_{\ell_0})})^*(N_{\ell_0 | \Pee^n})
= \Cal O_{\ell_0}(1)^{\oplus (n-1)}$ and that
$$N_{\ell_0 \times \Pee^{c_1+k-2r}|\Pee(\Cal O_{\ell_0}(1)^{\oplus k})} =
N_{\ell_0 \times \Pee^{c_1+k-2r}|\ell_0 \times \Pee^{k-1}}
= \Cal O_{\Pee^{c_1+k-2r}}(1)^{\oplus (2r-c_1-1)}.$$
Since $V|_{\ell_0} = \Cal O_{\ell_0}(1)^{\oplus k}
\oplus \oplus_{i=k+1}^r \Cal O_{\ell_0}(m_i)$,
$\xi|_{\ell_0 \times \Pee^{k-1}} =
\Cal O_{\ell_0}(1) \otimes \Cal O_{\Pee^{k-1}}(1)$ and
$$\align
N_{\Pee(\Cal O_{\ell_0}(1)^{\oplus k}) | \Pee(V|_{\ell_0})}
&= \oplus_{i=k+1}^r \Cal O_{\ell_0}(-m_i) \otimes
\xi|_{\ell_0 \times \Pee^{k-1}} \\
&= \oplus_{i=k+1}^r \Cal O_{\ell_0}(1-m_i) \otimes \Cal O_{\Pee^{k-1}}(1).\\
\endalign$$
Thus the previous two exact sequences are simplified to:
$$0 \to N_{\ell_0 \times \Pee^{c_1+k-2r}|\Pee(V|_{\ell_0})}
\to N_{\ell_0 \times \Pee^{c_1+k-2r}|\Pee(V)} \to
\Cal O_{\ell_0}(1)^{\oplus (n-1)} \to 0 \eqno (3.12)$$
$$0 \to \Cal O_{\Pee^{c_1+k-2r}}(1)^{\oplus (2r-c_1-1)}
\to N_{\ell_0 \times \Pee^{c_1+k-2r}|\Pee(V|_{\ell_0})} \to$$
$$\oplus_{i=k+1}^r \Cal O_{\ell_0}(1-m_i) \otimes
\Cal O_{\Pee^{c_1+k-2r}}(1) \to 0 \eqno (3.13)$$
Now (3.13) splits since for $k+1 \le i \le r$, we have $m_i \ge 2$ and
$$\align
&\quad \text{Ext}^1(\Cal O_{\ell_0}(1-m_i)
     \otimes \Cal O_{\Pee^{c_1+k-2r}}(1), \Cal O_{\Pee^{c_1+k-2r}}(1)) \\
&= H^1(\ell_0 \times \Pee^{c_1+k-2r}, \Cal O_{\ell_0}(m_i - 1)) = 0.\\
\endalign$$
Thus, the normal bundle $N_{\ell_0 \times \Pee^{c_1+k-2r}|\Pee(V|_{\ell_0})}$
is isomorphic to
$$\oplus_{i=k+1}^r \Cal O_{\ell_0}(1-m_i) \otimes
\Cal O_{\Pee^{c_1+k-2r}}(1) \oplus
\Cal O_{\Pee^{c_1+k-2r}}(1)^{\oplus (2r-c_1-1)},$$
and the exact sequence (3.12) becomes to the exact sequence:
$$0 \to \oplus_{i=k+1}^r \Cal O_{\ell_0}(1-m_i) \otimes
\Cal O_{\Pee^{c_1+k-2r}}(1) \oplus
\Cal O_{\Pee^{c_1+k-2r}}(1)^{\oplus (2r-c_1-1)} \to$$
$$N_{\ell_0 \times \Pee^{c_1+k-2r}|\Pee(V)} \to
\Cal O_{\ell_0}(1)^{\oplus (n-1)} \to 0 \eqno (3.14)$$
Restricting (3.14) to $\ell_0 \times \{ q \}$ and taking long exact
cohomology sequence result
$$\oplus_{i=k+1}^r H^1(\Cal O_{\ell_0}(1-m_i))
\otimes \Cal O_{\Pee^{c_1+k-2r}}(1)|_{q} \to$$
$$H^1((N_{\ell_0 \times \Pee^{c_1+k-2r}|\Pee(V)})|_{\ell_0 \times \{ q \}})
\to 0. \eqno (3.15)$$

Next, we determine $N_f$ and show that $h^1(N_f) \le c_1+k-2r$.
The two inclusions $\ell_0 \times \{ q \} \subset
\ell_0 \times \Pee^{c_1+k-2r} \subset \Pee(V)$ give an exact sequence
$$0 \to N_{\ell_0 \times \{ q \} | \ell_0 \times \Pee^{c_1+k-2r}}
\to N_{\ell_0 \times \{ q \}|\Pee(V)} \to
(N_{\ell_0 \times \Pee^{c_1+k-2r}|\Pee(V)})|_{\ell_0 \times \{ q \}} \to 0.$$
Since $N_{\ell_0 \times \{ q \} | \ell_0 \times
\Pee^{c_1+k-2r}} = T_{q, \Pee^{c_1+k-2r}}$, the above exact sequence becomes
$$0 \to T_{q, \Pee^{c_1+k-2r}} \to N_f \to
(N_{\ell_0 \times \Pee^{c_1+k-2r}|\Pee(V)})|_{\ell_0 \times \{ q \}} \to 0.
\eqno (3.16)$$
Thus, $h^1(N_f) =
h^1((N_{\ell_0 \times \Pee^{c_1+k-2r}|\Pee(V)})|_{\ell_0 \times \{ q \}})$.
By (3.15), we obtain
$$\align
h^1(N_f) &= h^1((N_{\ell_0 \times \Pee^{c_1+k-2r}|\Pee(V)})|_{\ell_0 \times
\{ q \}})  \le \sum_{i=k+1}^r h^1(\Cal O_{\ell_0}(1-m_i)) \\
&= \sum_{i=k+1}^r (m_i - 2) = c_1 + k - 2r.\\
\endalign$$

Finally, we show that $h^1(N_f) = c_1 + k - 2r$. It suffices to prove that
$h^1(N_f) \ge c_1 + k - 2r$. Since $\ell_0$ is a generic line in $\Pee^n$
and $V|_{\ell_0} = \Cal O_{\ell_0}(1)^{\oplus k}
\oplus \oplus_{i = k+1}^r \Cal O_{\ell_0}(m_i)$,
$\hbox{dim} \frak M(A, 0) = (2n + k)$ by Lemma 2.3 (ii).
Since $h^0(N_f)$ is the dimension of the Zariski tangent space of
$\frak M(A, 0)/PSL(2; \Cee)$ at $f$, $h^0(N_f) \ge (2n + k-3)$. Thus,
$$h^1(N_f) = h^0(N_f) - \chi(N_f) \ge (2n + k-3) - (2n + 2r-c_1-3)
= k + c_1 -2r.$$
Therefore, $h^1(N_f) = c_1 + k - 2r$. In particular, $h^1(N_f)$ is
independent of $f \in \frak M$. To obtain the obstruction bundle $COB$
over $\Pee^{c_1+k-2r}$, we notice that (3.15) gives
$$\oplus_{i=k+1}^r H^1(\Cal O_{\ell_0}(1-m_i)) \otimes
\Cal O_{\Pee^{c_1+k-2r}}(1)|_{q} \cong
H^1((N_{\ell_0 \times \Pee^{c_1+k-2r}|\Pee(V)})|_{\ell_0 \times \{ q \}}).$$
Thus by the exact sequence (3.16), we conclude that
$$\align
H^1(N_f)
&\cong H^1((N_{\ell_0 \times
  \Pee^{c_1+k-2r}|\Pee(V)})|_{\ell_0 \times \{ q \}}) \\
&\cong \oplus_{i=k+1}^r H^1(\Cal O_{\ell_0}(1-m_i))
  \otimes \Cal O_{\Pee^{c_1+k-2r}}(1)|_{q}. \tag 3.17
\endalign
$$
It follows that $COB = \Cal O_{\Pee^{c_1+k-2r}}(1)^{\oplus (c_1+k-2r)}$.
By (3.5), we obtain
$$a_0 = \Phi_{(A, 0)}(h, h_n, h_n \xi_{2r-c_1-1})
= e(COB) = 1. \qed$$

\section{4. Proof of Theorem B}

In this section, we prove Theorem B which we restate below.

\theorem{4.1} {\rm (i)} Let $V$ be a rank-$r$ ample bundle over $\Pee^n$.
Assume either $c_1 \le n$ or $c_1 \le (n + r)$ and
$V\otimes {\Cal O}_{\Pee^n}(-1)$ is nef so that $\Pee(V)$ is Fano.
Then the quantum cohomology $H^*_{\omega}(\Pee(V); \Zee)$ is
the ring generated by $h$ and $\xi$ with two relations
$$h^{n+1} = \sum_{i+j \le (c_1 - r)} a_{i,j} \cdot h^i \cdot \xi^j
\cdot e^{-t(n+1-i-j)} \eqno (4.2)$$
$$\sum_{i=0}^r (-1)^i c_i \cdot h^i \cdot \xi^{r-i} = e^{-tr}
+ \sum_{i+j \le (c_1-n-1)} b_{i,j} \cdot h^i \cdot \xi^{j}
\cdot e^{-t(r-i-j)} \eqno (4.3)$$
where the coefficients $a_{i, j}$ and $b_{i, j}$ are integers depending on $V$;

{\rm (ii)} If we further assume that $c_1<2r$,
then the leading coefficient $a_{0,c_1-r}=1$.
\endproclaim
\noindent
{\it Proof.} (i) First, we determine the first relation $f_\omega^1$
in (1.3). By Lemma 3.7,
$$h \cdot h_p = h_{p +1} + \sum_{A \in H_2'} (h \cdot h_p)_A \cdot
e^{t K_{\Pee(V)}(A)} \eqno (4.4)$$
where $p \ge 1$ and $H_2'$ stands for $H_2(\Pee(V); \Zee) -
\Zee \cdot (h_n\xi_{r-2})_*$. Thus,
$$h^{n - p} \cdot h_{p +1} = h^{n - p + 1} \cdot h_{p} - \sum_{A \in H_2'}
h^{n - p} \cdot (h \cdot h_p)_A \cdot e^{t K_{\Pee(V)}(A)}.$$
If $(h \cdot h_p)_A \ne 0$, then $A = [E]$ for some effective curve $E$.
So $a = h(A) \ge 0$. Since $A \in H_2'$, $a \ge 1$.
We claim that $-K_{\Pee(V)}(A) \ge (n+1-c_1 + r)$ with equality
if and only if $A = [h_{n-1}\xi_{r-1} + (1 - c_1) h_{n}\xi_{r-2}]_*
\overset \hbox{def} \to = A_2$. Indeed, if $c_1 \le n$,
then $-K_{\Pee(V)}(A) = (n+1-c_1)a + r \cdot \xi(A) \ge (n+1-c_1 + r)$
with equality if and only if $a = \xi(A) = 1$, that is,
if and only if $A = A_2$; if $c_1 \le (n + r)$ and $(\xi -h)$ is nef,
then again $-K_{\Pee(V)}(A) = (n+1+r-c_1)a + r \cdot (\xi - h)(A)
\ge (n+1-c_1 + r)$ with equality if and only if $a = 1$ and
$(\xi - h)(A) = 0$, that is, if and only if $A = A_2$.
Thus, $\hbox{deg}((h \cdot h_p)_A) = 1 + p + K_{\Pee(V)}(A)
\le (p - n +c_1 -r)$,
and $\hbox{deg}(h^{n - p} \cdot (h \cdot h_p)_A) \le (c_1 - r)$.
Using induction on $p$ and keeping track of the exponential
$e^{t K_{\Pee(V)}(A)}$, we obtain
$$0 = h_{n+1} = h^{n+1} - \sum_{i+j \le (c_1 - r)} a_{i,j} \cdot
h^i \cdot \xi^j \cdot e^{-t(n+1-i-j)}.$$
Therefore, the first relation $f_\omega^1$ for the quantum cohomology ring is:
$$h^{n+1} = \sum_{i+j \le (c_1 - r)} a_{i,j} \cdot h^i \cdot \xi^j \cdot
e^{-t(n+1-i-j)}.$$

Next, we determine the second relation $f_\omega^2$ in (1.3).
We need to compute the quantum product $h^i \cdot \xi^{r-i}$
for $0 \le i \le r$. First, we calculate the quantum product $\xi^r$.
Note that if $A = (bh_{n}\xi_{r-2})_*$ with $b \ge 1$,
then $-K_{\Pee(V)}(A) = br \ge r$ with
$-K_{\Pee(V)}(A) = r$ if and only if
$A = (h_{n}\xi_{r-2})_* \overset \hbox{def} \to = A_1$.
Thus for $p \ge 1$,
$$\xi \cdot \xi_p =
\cases
\xi_{p+1} + \sum_{A \in H_2'} (\xi \cdot \xi_p)_A \cdot e^{t K_{\Pee(V)}(A)},
  &\text{if $p < r -1$}\\
\xi_{r} + (\xi \cdot \xi_{r-1})_{A_1} \cdot e^{-tr} +
\sum_{A \in H_2'} (\xi \cdot \xi_{r-1})_A \cdot e^{t K_{\Pee(V)}(A)},
  &\text{if $p = r -1$.}\\
\endcases
$$
Note that $(\xi \cdot \xi_{r-1})_{A_1}$ is of
degree zero; by Lemma 3.6, we obtain $(\xi \cdot \xi_{r-1})_{A_1} =
\Phi_{(A_1, 0)}(\xi, \xi_{r-1}, h_{n}\xi_{r-1}) = 1$. Therefore for $p \ge 1$,
$$\xi \cdot \xi_p =
\cases
\xi_{p+1} + \sum_{A \in H_2'} (\xi \cdot \xi_p)_A \cdot e^{t K_{\Pee(V)}(A)},
  &\text{if $p < r -1$}\\
\xi_{r} + e^{-tr} +
\sum_{A \in H_2'} (\xi \cdot \xi_{r-1})_A \cdot e^{t K_{\Pee(V)}(A)},
  &\text{if $p = r -1$.}\\
\endcases
\eqno (4.5)$$
Now, for $i \ge 1$ and $j \ge 1$ with $i + j \le r$, we have
$$h_i \cdot \xi_j =
\cases
h_i \xi_{j} + \sum_{A \in H_2'} (h_i \cdot \xi_j)_A \cdot e^{t K_{\Pee(V)}(A)},
  &\text{if $i + j < r$}\\
h_i \xi_{j} + (h_i \xi_{j})_{A_1} \cdot e^{-tr} +
\sum_{A \in H_2'} (h_i \cdot \xi_j)_A \cdot e^{t K_{\Pee(V)}(A)},
  &\text{if $i + j = r$;}\\
\endcases
$$
when $i + j = r$, $(h_i \xi_{j})_{A_1}$ is of degree zero; by Lemma 3.7,
we have $(h_i \cdot \xi_j)_{A_1} =
\Phi_{(A_1, 0)}(h_i, \xi_{r-i}, h_{n}\xi_{r-1}) = 0$.
Therefore for $i \ge 1$ and $j \ge 1$ with $i + j \le r$,
$$h_i \cdot \xi_j = h_i \xi_{j} + \sum_{A \in H_2'}
(h_i \cdot \xi_j)_A \cdot e^{t K_{\Pee(V)}(A)}. \eqno (4.6)$$
 From the proof of the first relation $f_\omega^1$, we see that
if $\alpha$ and $\beta$ are homogeneous elements in
$H^*(\Pee(V); \Zee)$ with $\hbox{deg}(\alpha) + \hbox{deg}(\beta) = m \le r$,
then $\hbox{deg}((\alpha \cdot \beta)_A) \le m - (n+1-c_1 + r)$
for $A \in H_2'$. Thus if $\gamma$ is a homogeneous element in
$H^*(\Pee(V); \Zee)$ with $\hbox{deg}(\gamma) = r-m$,
then $\hbox{deg}(\gamma \cdot (\xi \cdot \xi_p)_A) \le (c_1 - n - 1)$.
Since $\sum_{i=0}^r (-1)^i c_i \cdot h_i \xi_{r-i} = 0$,
it follows from (4.4), (4.5), and (4.6) that
the second relation $f_\omega^2$ is
$$\sum_{i=0}^r (-1)^i c_i \cdot h^i \cdot \xi^{r-i} = e^{-tr}
+ \sum_{i+j \le (c_1-n-1)} b_{i,j} \cdot h^i \cdot \xi^{j}
\cdot e^{-t(r-i-j)}.$$

(ii) From the proof of the first relation in (i), we see that
$-K_{\Pee(V)}(A) \ge (n+1-c_1 + r)$ with equality
if and only if $A = A_2$;
moreover, the term $\xi^{c_1-r}$ can only come from the quantum correction
$(h \cdot h_n)_{A_2}$. Now
$$(h \cdot h_n)_{A_2} = (\sum_{i = 0}^{c_1 - r} a_i' h_i \xi_{c_1 - r - i})
\cdot e^{-t(n+1-c_1 + r)}$$
where $a_0' = \Phi_{(A_2, 0)}(h, h_n, h_n \xi_{2r-c_1-1})$.
Since $c_1 < 2r$, $(c_1 - r)< r$. By (4.4), (4.5), and (4.6), we conclude that
$h_i \xi_{c_1 - r - i} = h^i \cdot \xi^{c_1 - r - i} +
\text{(lower degree terms)}$.
Thus $a_{0,c_1-r} = a_0' = \Phi_{(A_2, 0)}(h, h_n, h_n \xi_{2r-c_1-1})$.
By Lemma 3.8, $a_{0,c_1-r} = 1.\qed$

It is understood that when $c_1 \le n$, then the summations
on the right-hand-sides of the second relations (4.3) and (4.9)
below do not exist.

Next, we shall sharpen the results in Theorem 4.1 by imposing additional
conditions on $V$. Let $V$ be a rank-$r$ ample vector bundle over $\Pee^n$.
Then $c_1 \ge r$. Thus if $c_1<2r$ and if
either $2c_1 \le (n + r)$ or $2c_1 \le (n + 2r)$ and
$V\otimes {\Cal O}_{\Pee^n}(-1)$ is nef, then the conditions in Theorem 4.1
are satisfied.

\corollary{4.7} {\rm (i)} Let $V$ be a rank-$r$ ample vector bundle over
$\Pee^n$ with $c_1<2r$. Assume that either $2c_1 \le (n + r)$ or
$2c_1 \le (n + 2r)$ and $V\otimes {\Cal O}_{\Pee^n}(-1)$ is nef so that
$\Pee(V)$ is a Fano variety. Then the first relation (4.2) is
$$h^{n+1}= \left ( \sum_{i=0}^{c_1-r} a_i \cdot h^i \cdot \xi^{c_1-r-i}
\right ) \cdot e^{-t(n+1+r-c_1)} \eqno (4.8)$$
where the integers $a_i$ depend on $V$. Moreover, $a_0 = 1$.

{\rm (ii)} Let $V$ be a rank-$r$ ample vector bundle over $\Pee^n$.
Assume that $2c_1 \le (2n + r+1)$ and $V\otimes {\Cal O}_{\Pee^n}(-1)$ is nef
so that $\Pee(V)$ is Fano. Then the second relation (4.3) is
$$\sum_{i=0}^r (-1)^i c_i \cdot h^i \cdot \xi^{r-i} = e^{-tr} +
\sum_{i=0}^{c_1-n-1} b_{i} \cdot h^i \cdot \xi^{c_1-n-1-i}
\cdot e^{-t(n+1+r-c_1)} \eqno (4.9)$$
where the integers $b_i$ depend on $V$.
\endproclaim
\proof (i) From the proof of Theorem 4.1 (i),
we notice that it suffices to show that the only homology class
$A \in H_2' = H_2(\Pee(V); \Zee) - \Zee \cdot (h_n\xi_{r-2})_*$
which has nonzero contributions to the quantum corrections in (4.4) is
$A = [h_{n-1}\xi_{r-1} + (1 - c_1) h_{n}\xi_{r-2}]_*
\overset \hbox{def} \to = A_2$. In other words,
if $A = (ah_{n-1}\xi_{r-1}+ b h_{n}\xi_{r-2})_*$ with $a \ne 0$
and if $\Phi_{(A, 0)}(h, h_p, \alpha) \ne 0$ for $1 \le p \le n$ and
$\alpha \in H^*(\Pee(V); \Zee)$, then $A = A_2$. First of all,
we show that $a = 1$. Suppose $a \ne 1$. Then $a \ge 2$. By (3.2),
$$\align
{1 \over 2} \cdot \hbox{deg}(\alpha)
&= (n+r-1) - K_{\Pee(V)}(A) - 1 - p\cr
&= (n+r-1) + [(n + 1-c_1)a + r \cdot \xi(A)] - 1-p\cr
&\ge \hbox{dim}(\Pee(V)) + [(n + 1-c_1)a + r \cdot \xi(A)] - 1-n.\cr
\endalign$$
If $2c_1 \le (n + r)$, then $c_1 \le n$,
and $[(n + 1-c_1)a + r \cdot \xi(A)] - 1-n \ge 2(n + 1-c_1) + r - 1-n > 0$.
If $2c_1 \le (n + 2r)$ and $(\xi - h)$ is nef,
then $c_1 \le n + r$, and $[(n + 1-c_1)a + r \cdot \xi(A)] - 1-n =
[(n + 1 + r-c_1)a + r \cdot (\xi-h)(A)] - 1-n \ge 2(n + 1 + r-c_1)-1 -n > 0$.
Thus, $[(n + 1-c_1)a + r \cdot \xi(A)] - 1-n > 0$,
and so $\hbox{deg}(\alpha)/2 > \hbox{dim}(\Pee(V))$. But this is absurd.
Next, we prove that $b = (1 - c_1)$, or equivalently, $\xi(A) = 1$.
Suppose $\xi(A) \ne 1$. Then $\xi(A) \ge 2$. By (3.2),
$$\align
{1 \over 2} \cdot \hbox{deg}(\alpha)
&= (n+r-1) + [(n + 1-c_1) + r \cdot \xi(A)] - 1-p\cr
&\ge \hbox{dim}(\Pee(V)) + [(n + 1-c_1) + 2r] - 1-n\cr
&> \hbox{dim}(\Pee(V))\cr
\endalign$$
since $c_1 < 2r$. But once again this is absurd.

(ii) We follow the previous arguments for (i). Again it suffices to show that
if $A = (ah_{n-1}\xi_{r-1}+ b h_{n}\xi_{r-2})_*$ with $a \ne 0$
and if $\Phi_{(A, 0)}(\alpha_1, \alpha_2, \alpha) \ne 0$ for
some $\alpha_1, \alpha_2, \alpha \in H^*(\Pee(V); \Zee)$
with $\hbox{deg}(\alpha_1) + \hbox{deg}(\alpha_2) \le r$, then $A = A_2$.
Indeed, if $a \ne 1$ or if $a = 1$ but $\xi(A) \ne 1$,
then we must have $\hbox{deg}(\alpha)/2 > \hbox{dim}(\Pee(V))$.
But this is impossible. Therefore, $a = 1$ and $\xi(A) = 1$. So $A = A_2$.
\endproof

Now we discuss the relation between the quantum corrections and
the extremal rays of the Fano variety $\Pee(V)$.
Let $V$ be a rank-r ample vector bundle over $\Pee^n$
with $c_1<2r$ and $2c_1 \le (n + r)$. By (4.8) and (4.3),
the quantum cohomology ring $H^*_{\omega}(\Pee(V); \Zee)$ is
the ring generated by $h$ and $\xi$ with two relations
$$h^{n+1}= \left ( \sum_{i=0}^{c_1-r} a_i \cdot h^i \cdot \xi^{c_1-r-i}
\right ) \cdot e^{-t(n+1+r-c_1)} \eqno (4.10)$$
$$\sum_{i=0}^r (-1)^i c_i \cdot h^i \cdot \xi^{r-i} = e^{-tr}.
\eqno (4.11)$$
 From the proof of Theorem 4.1 (i), we notice that the quantum correction to
the second relation (4.11) comes from the homology class
$A_1 = (h_n\xi_{r-2})_*$ which is represented by the lines in
the fibers of $\pi: \Pee(V) \to \Pee^n$. Also,
we notice from the proof of Corollary 4.7 (i) that the quantum correction
to the first relation (4.10) comes from the homology class
$A_2 = [h_{n-1}\xi_{r-1} + (1 - c_1) h_{n}\xi_{r-2}]_*$;
from the proof of Lemma 3.8, $A_2$ can be represented by
a smooth rational curve isomorphic to lines in $\Pee^n$ via $\pi$.
Now $A_1$ generates one of the two extremal rays of
$\Pee(V)$. It is unclear whether $A_2$ generates the other extremal ray.
By Lemma 2.2 (iii), if we further assume that $2c_1 \le (n + 1)$,
then indeed $A_2$ generates the other extremal ray of $\Pee(V)$.
By Lemma 2.2 (ii), $A_2$ generates the other extremal ray of
$\Pee(V)$ if and only if $(\xi - h)$ is nef, that is,
$V \otimes \Cal O_{\Pee^n}(-1)$ is a nef vector bundle over $\Pee^n$.

\section{5. Direct sum of line bundles over $\Pee^n$}

In this section, we partially verify Batyrev's conjecture on
the quantum cohomology of projective bundles associated to
direct sum of line bundles over $\Pee^n$.
We shall use (3.5) to compute the necessary Gromov-Witten invariants.
Our first step is to recall some standard materials for
the Grassmannian $G(2, n+1)$ from \cite{3}.
Then we determine certain obstruction bundle and its Euler class.
Finally we proceed to determine the first and second relations
for the quantum cohomology.

On the Grassmannian $G(2, n+1)$, there exists a tautological exact sequence
$$0 \to S \to (\Cal O_{G(2, n+1)})^{\oplus (n + 1)} \to Q \to 0 \eqno (5.1)$$
where the sub- and quotient bundles $S$ and $Q$ are of
rank $2$ and $(n -1)$ respectively. Let $\alpha$ and $\beta$
be the virtual classes such that $\alpha + \beta = -c_1(S)$
and $\alpha \beta = c_2(S)$. Then
$$\hbox{cl}(\{ \ell \in G(2, n+1)| \ell \cap h_p \ne \emptyset \})
= {\alpha^p - \beta^p \over \alpha - \beta} \eqno (5.2)$$
where $\hbox{cl}(\cdot)$ denotes the fundamental class and
$h_p$ stands for a fixed linear subspace of $\Pee^n$ of codimention $p$.
If $P(\alpha, \beta)$ is a symmetric homogeneous polynomial of
degree $(2n - 2)$ (so that $P(\alpha, \beta)$ can be written as
a polynomial of maximal degree in the Chern classes of the bundle $S$),
then we have
$$\int_{G(2, n+1)} P(\alpha, \beta) =
\left ( \text{the coefficient of } \alpha^n \beta^n \text{ in}
-{1 \over 2}(\alpha - \beta)^2 P(\alpha, \beta) \right ).  \eqno (5.3)$$
Let $F_n = \{ (x, \ell) \in \Pee^n \times G(2, n+1)| x \in \ell \}$,
and $\pi_1$ and $\pi_2$ are the two natural projections from
$F_n$ to $\Pee^n$ and $G(2, n+1)$ respectively.
Then $F_n = \Pee(S^*)$ where $S^*$ is the dual bundle of $S$,
and $(\pi_1^*\Cal O_{\Pee^n}(1))|_{F_n}$ is
the tautological line bundle over $F_n$.
Let $\hbox{Sym}^m(S^*)$ be the $m$-th symmetric product of $S^*$.
Then for $m \ge 0$,
$$\pi_{2*} (\pi_1^*\Cal O_{\Pee^n}(m)|_{F_n})
\cong \hbox{Sym}^m(S^*).  \eqno (5.4)$$
By the duality theorem for higher direct image sheaves
(see p.253 in \cite{7}),
$$\align
R^1 \pi_{2*} (\pi_1^*\Cal O_{\Pee^n}(-m)|_{F_n})
&\cong (\pi_{2*} (\pi_1^*\Cal O_{\Pee^n}(m - 2)|_{F_n}))^*
  \otimes (\hbox{det}S^*)^* \\
&\cong \hbox{Sym}^{m-2}(S) \otimes (\hbox{det}S) \tag 5.5 \\
\endalign$$

Now, let $V = \oplus_{i = 1}^r \Cal O_{\Pee^n}(m_i)$ where
$1 = m_1 = \ldots = m_k < m_{k+1} \le \ldots \le m_r$.
Assume that $k \ge 1$ and $\Pee(V)$ is Fano.
Then the two extremal rays of $\Pee(V)$ are generated by the two classes
$A_1 = (h_{n}\xi_{r-2})_*$ and
$A_2 = [h_{n-1}\xi_{r-1}+ (1 - c_1) h_{n}\xi_{r-2}]_*$. From the proof
of Lemma 2.3 (ii), we see that
$$\frak M(A_2, 0)/PSL(2; \Cee) = G(2, n+1) \times \Pee^{k -1}.  \eqno (5.6)$$
Let a morphism $f \in \frak M(A_2, 0)$ be induced by some surjective map
$V|_\ell \to \Cal O_\ell(1) \to 0$ such that
the image $\hbox{Im}(f)$ of $f$ is of the form
$$\hbox{Im}(f) = \ell \times \{q\} \subset \ell \times \Pee^{k-1}
\subset \Pee^n \times \Pee^{k-1}. $$
Then by arguments similar to the proof of (3.17), we have
$$H^1(N_f) \cong \oplus_{u=k+1}^r H^1(\Cal O_{\ell}(1-m_u))
\otimes \Cal O_{\Pee^{k-1}}(1)|_{q}. \eqno (5.7)$$
It follows that the obstruction bundle $COB$ over
$\frak M(A_2, 0)/PSL(2; \Cee)$ is
$$COB \cong \oplus_{u=k+1}^r R^1 \pi_{2*} (\pi_1^*\Cal
O_{\Pee^n}(1-m_u)|_{F_n})
\otimes \Cal O_{\Pee^{k-1}}(1). \eqno (5.8)$$
Since $c_1(S) = -(\alpha + \beta)$ and $c_2(S) = \alpha \beta$,
we obtain from (5.5) the following.

\lemma{5.9} The Euler class of the obstruction bundle $COB$ is
$$e(COB) = \prod_{u = k +1}^r \prod_{v=0}^{m_u - 3}
[(1+v)(-\alpha) + (m_u - 2 -v)(-\beta) + \tilde h] \eqno (5.10)$$
where $\tilde h$ stands for the hyperplane class in $\Pee^{k-1}$. \qed
\endproclaim

Next assuming $c_1 < 2r$, we shall compute the Gromov-Witten invariant
$$W_i \overset \hbox{def} \to = \Phi_{(A_2, 0)}(h_{\tilde n},
h_{n+1 - \tilde n}, h_{n-i}\xi_{2r-c_1 - 1+i}) \eqno (5.11)$$
where $0 \le i \le (c_1 - r)$ and $\tilde n = \left [{n+1 \over 2} \right ]$
is the largest integer $\le (n+1)/2$.

\lemma{5.12} Assume $c_1 < \text{min}(2r, (n+1+2r)/2)$ and
$0 \le i \le (c_1 - r)$. Then $W_i$ is the coefficient of $t^i$
in the power series expansion of
$$\prod_{u = 1}^r (1 - m_u t)^{m_u-2}.$$
\endproclaim
\noindent
{\it Proof.} Note that the restriction of $\xi$ to $\Pee^n \times \Pee^{k-1}
= \Pee(\Cal O_{\Pee^n}(1)^{\oplus k})$ is $(h + \tilde h)$. Thus,
$$\align
h_{n-i}\xi_{2r-c_1 - 1+i}|_{\Pee^n \times \Pee^{k-1}}
&= \sum_{j = 0}^{2r-c_1-1+i} {2r-c_1 - 1+i \choose j}
h_{n-i+j} \tilde h_{2r-c_1-1+i-j}\\
&= \sum_{j = 0}^{i} {2r-c_1 - 1+i \choose j}
h_{n-i+j} \tilde h_{2r-c_1-1+i-j}. \\
\endalign$$
So by (3.5) (replacing $\frak M(A_2, 0)$ by
$\frak M(A_2, 0)/PSL(2; \Cee)$), (5.2), and Lemma 5.9,
$$W_i = \int_{G(2, n+1) \times \Pee^{k-1}} \tilde P(\alpha, \beta)
\eqno (5.13)$$
where $\tilde P(\alpha, \beta)$ is the symmetric homogeneous polynomial of
degree $(2n - 2) + (k-1)$:
$$\align
\tilde P(\alpha, \beta) &= {\alpha^{\tilde n} - \beta^{\tilde n}
\over \alpha - \beta} \cdot {\alpha^{n + 1- \tilde n} -
\beta^{n + 1- \tilde n} \over \alpha - \beta} \\
&\qquad \cdot \sum_{j = 0}^{i} {2r-c_1 - 1+i \choose j}
{\alpha^{n-i+j} - \beta^{n-i+j} \over \alpha - \beta}
\cdot \tilde h_{2r-c_1-1+i-j} \\
&\qquad \cdot \prod_{u = k +1}^r \prod_{v=0}^{m_u - 3}
[(1+v)(-\alpha) + (m_u - 2 -v)(-\beta) + \tilde h] \\
&= \sum_{j = 0}^{i} {2r-c_1 - 1+i \choose j}
{\alpha^{n+1} - \alpha^{n + 1- \tilde n}\beta^{\tilde n} - \alpha^{\tilde n}
\beta^{n + 1- \tilde n} + \beta^{n+1} \over (\alpha - \beta)^2}\\
&\qquad \cdot \sum_{t= 0}^{n-i+j-1} \alpha^{t} \beta^{n-i+j-1-t} \cdot
\tilde h_{2r-c_1-1+i-j}\\
&\qquad \cdot \prod_{u = k +1}^r \prod_{v=0}^{m_u - 3}
[(1+v)(-\alpha) + (m_u - 2 -v)(-\beta) + \tilde h]. \\
\endalign$$
By (5.3) and (5.13), we conclude from straightforward manipulations that:
$$\align
W_i &= \sum_{j = 0}^{i} {2r-c_1 - 1+i \choose j} \cdot (-1)^{i-j} \\
&\qquad \cdot \sum_{j_{k+1} + \ldots + j_r = i -j} \quad
\prod_{u = k+1}^r {m_u-2 \choose j_u} (m_u -1)^{j_u} \\
 &= \sum_{j = 0}^{i} {2r-c_1 - 1+i \choose i-j} \cdot (-1)^{j} \\
&\qquad \cdot \sum_{j_{k+1} + \ldots + j_r = j} \quad
\prod_{u = k+1}^r {m_u-2 \choose j_u} (m_u -1)^{j_u}. \\
\endalign$$
Thus $W_i$ is the coefficient of $t^i$ in the polynomial
$$\align
&\qquad (1+t)^{2r-c_1 - 1+i} \cdot \prod_{u = k+1}^r [1 - (m_u -1) t]^{m_u-2}\\
&= (1+t)^{2r-c_1 - 1+i} \cdot \prod_{u = k+1}^r [(1+t) - m_u t]^{m_u-2}\\
&= (1+t)^{2r-c_1 - 1+i} \cdot \sum_{j=0}^{c_1-2r+k}
\sum_{j_{k+1} + \ldots + j_r = j} \\
&\qquad \cdot \prod_{u = k+1}^r {m_u-2 \choose j_u} (-m_ut)^{j_u}
\cdot (1+t)^{m_u-2-j_u}\\
&= \sum_{j=0}^{c_1-2r+k} \sum_{j_{k+1} + \ldots + j_r = j} \quad
\prod_{u = k+1}^r {m_u-2 \choose j_u} (-m_ut)^{j_u} \cdot (1+t)^{i+k-1-j}\\
\endalign$$
since $\sum_{u = k+1}^r (m_u - 2 - j_u) = c_1 -2r+k - j$.
So $W_i$ is the coefficient of $t^i$ in
$$\align
\prod_{u = k+1}^r (1-m_u t)^{m_u-2} \cdot \sum_{j=0}^{+ \infty}
{j +k-1 \choose k-1} t^j
&= \prod_{u = k+1}^r (1-m_u t)^{m_u-2} \cdot {1 \over (1 - t)^k}\\
&= \prod_{u = 1}^r (1 - m_u t)^{m_u-2}. \qed\\
\endalign$$

\proposition{5.14} Let $V = \oplus_{i = 1}^r \Cal O_{\Pee^n}(m_i)$ where
$m_i \ge 1$ for each $i$ and
$$\sum_{i=1}^r m_i < \text{min}(2r, (n+1+2r)/2).$$
Then the first relation $f^1_\omega$ for the quantum cohomology ring
$H^*_{\omega}(\Pee(V); \Zee)$ is
$$h^{n+1} = \prod_{u = 1}^r (\xi - m_u h)^{m_u-1} \cdot
e^{-t(n+1+r-\sum_{i=1}^r m_i)}. \eqno (5.15)$$
\endproclaim
\noindent
{\it Proof.} We may assume that
$1 = m_1 = \ldots = m_k < m_{k+1} \le \ldots \le m_r$.
Since the conclusion clearly holds when $k = r$, we also assume that $k < r$.
Let $c_1 = \sum_{i=1}^r m_i$.
Notice that the conditions in Corollary 4.7 (i) are satisfied. Thus,
$$h^{n+1}= \left ( \sum_{i=0}^{c_1-r} a_i \cdot h^i \cdot \xi^{c_1-r-i}
\right ) \cdot e^{-t(n+1+r-c_1)}.$$
More directly, putting $\tilde n = \left [{n+1 \over 2} \right ]$,
then $\tilde n < -K_{\Pee(V)}(A_2) = (n+1+r -c_1)$,
and $(n+1-\tilde n) < -K_{\Pee(V)}(A_2)$
unless $n$ is even and $c_1 = (n + 2r)/2$. From the proofs in Theorem 4.1
and Corollary 4.7 (i) for the first relation $f^1_\omega$,
we have $h^{\tilde n} = h_{\tilde n}$,
and $h^{n+1-\tilde n} = h_{n+1-\tilde n}$ unless $n$ is even and
$c_1 = (n + 2r)/2$. Moreover, if $n$ is even and $c_1 = (n + 2r)/2$,
then $h^{n+1-\tilde n} = h \cdot h^{n-\tilde n} = h \cdot h_{n-\tilde n}
= h_{n+1-\tilde n} + (h \cdot h_{n-\tilde n})_{A_2} \cdot e^{-t(n+1+r-c_1)}$.
Since $(h \cdot h_{n-\tilde n})_{A_2}$ is of degree zero,
$(h \cdot h_{n-\tilde n})_{A_2} =
\Phi_{(A_2, 0)}(h, h_{n-\tilde n}, h_n\xi_{r-1})$. Since $1 \le k <r$,
we can choose a point $q_0$ in $\Pee(V)$ representing the homology class
$(h_n\xi_{r-1})_*$ such that the point $q_0$ is not contained in
the $(k-1)$-dimensional linear subspace
$$\Pee^{k -1} = \Pee((\Cal O_{\Pee^n}(1)^{\oplus k})|_{\pi(q_0)})
\subset \Pee(V|_{\pi(q_0)}) \cong \Pee^{r-1}.$$
Note that for every $f \in \frak M(A_2, 0)$,
$\hbox{Im}(f) = \ell \times \{ q \}$ for some line $\ell \subset \Pee^n$
and some point $q \in \Pee^{k -1}$. Thus $\hbox{Im}(f)$ can not pass $q_0$.
As in the proof of Lemma 3.7, we conclude that
$\Phi_{(A_2, 0)}(h, h_{n-\tilde n}, h_n\xi_{r-1}) = 0$.
Therefore, $h^{n+1-\tilde n} = h_{n+1-\tilde n}$. So
$$h^{n+1} = h^{\tilde n} \cdot h^{n+1-\tilde n}
= h_{\tilde n} \cdot h_{n+1-\tilde n}.$$
By similar arguments in the proofs of Theorem 4.1 and Corollary 4.7 (i)
for the first relation $f^1_\omega$,
we see that if $(h_{\tilde n} \cdot h_{n+1-\tilde n})_A \ne 0$,
then $A = 0, A_2$. Thus
$$h^{n+1} = h_{n+1} + (h_{\tilde n} \cdot h_{n+1-\tilde n})_{A_2} \cdot
e^{-t(n+1+r-c_1)} = (h_{\tilde n} \cdot h_{n+1-\tilde n})_{A_2}
\cdot e^{-t(n+1+r-c_1)}.$$
So it suffices to show that $(h_{\tilde n} \cdot h_{n+1-\tilde n})_{A_2}
= \prod_{u = 1}^r (\xi - m_u h)^{m_u-1}$. Note that
$$\prod_{u = 1}^r (\xi - m_u h)^{m_u-1} =
\prod_{u = 1}^r (\xi - m_u h)_{m_u-1}$$
where the right-hand-side stands for the product in the ordinary cohomology.
Thus we need to show that $(h_{\tilde n} \cdot h_{n+1-\tilde n})_{A_2}
= \prod_{u = 1}^r (\xi - m_u h)_{m_u-1}$, or equivalently,
$$\Phi_{(A_2, 0)}(h_{\tilde n}, h_{n+1 - \tilde n}, h_{n-i}\xi_{2r-c_1 - 1+i})
= \prod_{u = 1}^r (\xi - m_u h)_{m_u-1} h_{n-i}\xi_{2r-c_1 - 1+i}
\eqno (5.16)$$
for $0 \le i \le (c_1 - i)$. The left-hand-side of (5.16) is computed
in Lemma 5.12.

Denote the right-hand-side of (5.16) by $\tilde W_i$.
Let $s_i$ be the $i$-th Segre class of $V$. Then
we have $s_i = (-1)^i \cdot \sum_{j_1 + \ldots + j_r = i}
\prod_{u=1}^r m_u^{j_u}$ and
$$\sum_{i=0}^{+ \infty} (-1)^{i} s_{i} t^{i} =
\prod_{u=1}^r {1 \over 1-m_u t}. \eqno (5.17)$$
Moreover from the second relation in (1.1), we obtain for $i \ge r$,
$$\xi_i = (-1)^{i - (r-1)} s_{i-(r-1)} \xi_{r-1} +
\text{ (terms with exponentials of } \xi \text{ less than } (r-1)).$$
It follows from the right-hand-side of (5.16) that
$\tilde W_i$ is equal to
$$\align
&\quad \sum_{j=0}^{c_1-r} \quad \sum_{j_1 + \ldots + j_r = j} \quad
\prod_{u=1}^r {m_u-1 \choose j_u} \xi_{m_u-1-j_u} (-m_uh)_{j_u}
h_{n-i}\xi_{2r-c_1 - 1+i} \\
&= \sum_{j=0}^{i} \quad \sum_{j_1 + \ldots + j_r = j} \quad \prod_{u=1}^r
{m_u-1 \choose j_u} (-m_u)^{j_u} h_{n-i+j}\xi_{r - 1+i-j} \\
&= \sum_{j=0}^{i} (-1)^{i-j} s_{i-j} \quad \sum_{j_1 + \ldots + j_r = j}
\quad \prod_{u=1}^r {m_u-1 \choose j_u} (-m_u)^{j_u}. \\
\endalign$$
Therefore, the formal power series $\sum_{i=0}^{+ \infty} \tilde W_i t^i$
is equal to
$$\align
&\quad \sum_{i=0}^{+ \infty} \sum_{j=0}^{i} (-1)^{i-j} s_{i-j} t^{i-j}
\quad \sum_{j_1 + \ldots + j_r = j} \quad
\prod_{u=1}^r {m_u-1 \choose j_u} (-m_u t)^{j_u} \\
&= \sum_{j=0}^{+ \infty} \sum_{i=j}^{+ \infty} (-1)^{i-j} s_{i-j} t^{i-j}
\quad \sum_{j_1 + \ldots + j_r = j} \quad
\prod_{u=1}^r {m_u-1 \choose j_u} (-m_u t)^{j_u} \\
&= \sum_{j=0}^{+ \infty} \sum_{i=0}^{+ \infty} (-1)^{i} s_{i} t^{i}
\quad \sum_{j_1 + \ldots + j_r = j} \quad
\prod_{u=1}^r {m_u-1 \choose j_u} (-m_u t)^{j_u} \\
&= \sum_{j=0}^{+ \infty} \prod_{u=1}^r {1 \over 1-m_u t}
\quad \sum_{j_1 + \ldots + j_r = j} \quad
\prod_{u=1}^r {m_u-1 \choose j_u} (-m_u t)^{j_u} \\
&= \prod_{u=1}^r {1 \over 1-m_u t} \sum_{j=0}^{+ \infty}
\quad \sum_{j_1 + \ldots + j_r = j} \quad
\prod_{u=1}^r {m_u-1 \choose j_u} (-m_u t)^{j_u} \\
&= \prod_{u=1}^r {1 \over 1-m_u t} \quad \prod_{u=1}^r (1-m_u t)^{m_u-1}\\
&= \prod_{u=1}^r (1-m_u t)^{m_u-2}\\
\endalign$$
where we have applied (5.17) in the third equality.
By Lemma 5.12, $\tilde W_i = W_i$ for $0 \le i \le (c_1 - r)$.
Hence the formule (5.16) and (5.15) hold.
\qed

It turns out that under certain conditions on the integers $m_i$,
the second relation $f^2_\omega$ for the quantum cohomology ring
$H^*_{\omega}(\Pee(V); \Zee)$ is much easier to be determined.
Note that the second relation $f^2$ in (1.1) can be rewritten as
$$\prod_{i=1}^r (\xi-m_i h) = 0  \eqno (5.18)$$
where the left-hand-side stands for the product in the ordinary cohomology
ring.

\proposition{5.19} Let $V = \oplus_{i = 1}^r \Cal O_{\Pee^n}(m_i)$ where
$m_i \ge 1$ for each $i$, $m_i = 1$ for some $i$, and
$\sum_{i=1}^r m_i < (2n+2+r)/2$.
Then the second relation $f^2_\omega$ for the quantum cohomology ring
$H^*_{\omega}(\Pee(V); \Zee)$ is
$$\prod_{i=1}^r (\xi-m_i h) = e^{-tr} \eqno (5.20)$$
where the left-hand-side stands for the product in the quantum cohomology ring.
\endproclaim
\noindent
{\it Proof.} We may assume that
$1 = m_1 = \ldots = m_k < m_{k+1} \le \ldots \le m_r$. So $k \ge 1$.
We notice that the conditions in Corollary 4.7 (ii)
are satisfied. From the proofs of Theorem 4.1 (i) and Corollary 4.7 (ii),
we see that the quantum corrections to the second relation (5.18)
can only come from the classes $A_1, A_2$; moreover,
the quantum correction from $A_1$ is $e^{-tr}$.
Thus it suffices to show that the quantum correction from $A_2$
is zero. In view of (3.3), it suffices to show that
$$\Phi_{(A_2, 0)}(\xi-m_1 h, \ldots, \xi-m_r h, \alpha) = 0$$
for every $\alpha \in H^*(\Pee(V); \Zee)$. For $1 \le i \le r$,
let $V_i$ be the subbundle of $V$:
$$V_i = \Cal O_{\Pee^n}(m_1) \oplus \ldots \oplus
\Cal O_{\Pee^n}(m_{i-1}) \oplus \Cal O_{\Pee^n}(m_{i+1})
\oplus \ldots \oplus \Cal O_{\Pee^n}(m_r),$$
and let $B_i = \Pee(V_i)$ be
the codimension-$1$ subvariety of $\Pee(V)$ induced by the projection
$V \to V_i \to 0$. Then the fundamental class of $B_i$ is $(\xi-m_i h)$.
As in the proof of Lemma 3.7, we need only to show that if
$f \in \frak M(A_2, 0)$, then the image $\hbox{Im}(f)$ can not intersect
with $B_1, \ldots, B_r$ simultaneously. In fact, we will show that
$\hbox{Im}(f)$ can not intersect with $B_1, \ldots, B_k$ simultaneously.
Indeed, $\hbox{Im}(f)$ is of the form
$$\hbox{Im}(f) = \ell \times \{q\} \subset \ell \times \Pee^{k-1}
\subset \Pee^n \times \Pee^{k-1} = \Pee(\Cal O_{\Pee^n}(1)^{\oplus k})$$
for some line $\ell \subset \Pee^n$,
and $B_i|_{\pi^{-1}(\ell)} = \Pee(V_i|_\ell)$.
Put $p = \pi(q) \in \Pee^n$, and
$$V|_p = \oplus_{i=1}^k \Cee \cdot e_i \oplus
(\oplus_{i=k+1}^r \Cal O_{\Pee^n}(m_i)|_p)$$
where $e_i$ is a global section of
$\Cal O_{\Pee^n}(m_i) = \Cal O_{\Pee^n}(1)$ for $i \le k$.
Now the point $q$ is identified with $\Cee \cdot v$ for some nonzero vector
$v \in \oplus_{i=1}^k \Cee \cdot e_i$. Let $v = \sum_{i=1}^k a_i e_i$.
Since $\ell \times \{q\}$ and $B_i$ ($1 \le i \le k$) intersect,
the one-dimensional vector space $\Cee \cdot v$ is also contained
in $(V_i)|_p$. It follows that $a_i = 0$ for every $i$ with
$1 \le i \le k$. But this is impossible since $v$ is a nonzero vector.
\qed

In summary, we partially verify Batyrev's conjecture.

\theorem{5.21} Let $V = \oplus_{i = 1}^r \Cal O_{\Pee^n}(m_i)$ where
$m_i \ge 1$ for each $i$ and
$$\sum_{i=1}^r m_i < \text{min}(2r, (n+1+2r)/2, (2n+2+r)/2).$$
Then the quantum cohomology $H^*_{\omega}(\Pee(V); \Zee)$ is
generated by $h$ and $\xi$ with relations
$$h^{n+1}=\prod^r_{i=1}(\xi-m_ih)^{m_i-1} \cdot
e^{-t(n+1+r-\sum_{i=1}^r m_i)} \qquad {and} \qquad
\prod^r_{i=1}(\xi-m_ih) = e^{-tr}.$$
\endproclaim
\proof Follows immediately from Propositions 5.14 and 5.19.
\endproof

\section{6. Examples}

In this section, we shall determine the quantum cohomology of $\Pee(V)$
for ample bundles $V$ over $\Pee^n$ with $2 \le r \le n$ and $c_1 = r + 1$.
In these cases, $V|_\ell = \Cal O_\ell(1)^{\oplus (r - 1)} \oplus
\Cal O_\ell(2)$ for every line $\ell \subset \Pee^n$. In particular,
$V$ is a uniform bundle. If $r < n$, then by the Theorem 3.2.3 in \cite{10},
$V = \Cal O_{\Pee^n}(1)^{\oplus (r - 1)} \oplus \Cal O_{\Pee^n}(2)$;
if $r = n$, then by the results on pp.71-72 in \cite{10},
$V = \Cal O_{\Pee^n}(1)^{\oplus (n - 1)} \oplus \Cal O_{\Pee^n}(2)$
or $V = T_{\Pee^n}$ the tangent bundle of $\Pee^n$.
When $V = \Cal O_{\Pee^n}(1)^{\oplus (r - 1)} \oplus \Cal O_{\Pee^n}(2)$
with $r \le n$, the conditions in Theorem 5.21 are satisfied,
so the quantum cohomology ring $H^*_{\omega}(\Pee(V); \Zee)$ is
the ring generated by $h$ and $\xi$ with two relations
$$h^{n+1}= (\xi-2h) \cdot e^{-t(n+1+r-c_1)}
\qquad \hbox{and} \qquad (\xi-h)^{r-1} (\xi-2h) = e^{-tr}.$$

In the rest of this section, we compute the quantum cohomology of
$\Pee(T_{\Pee^n})$. It is well-known that $(\xi-h)$ is a nef divisor
on $\Pee(T_{\Pee^n})$, and the two extremal rays of $\Pee(T_{\Pee^n})$
are generated by $A_1 = (h_{n}\xi_{n-2})_*$ and
$A_2 = (h_{n-1}\xi_{n-1} - n h_{n}\xi_{n-2})_*$.
Moreover, $A_2$ is represented by smooth rational curves in $\Pee(T_{\Pee^n})$
induced by the surjective maps $T_{\Pee^n}|_\ell \to \Cal O_\ell(1) \to 0$
for lines $\ell \subset \Pee^n$. Since $c_1 = n+1$ and $n \ge 2$,
the assumptions in Corollary 4.7 are satisfied,
so the quantum cohomology ring $H^*_{\omega}(\Pee(T_{\Pee^n}); \Zee)$
is the ring generated by $h$ and $\xi$ with two relations
$$h^{n+1}= (a_1 h + \xi) \cdot e^{-tn} \quad \text{and} \quad \sum_{i=0}^n
(-1)^i c_i \cdot h^i \cdot \xi^{n-i} = (1 + b_0) \cdot e^{-tn}. \eqno (6.1)$$
More precisely, putting $H_2' = H_2(\Pee(V); \Zee) -
\Zee \cdot (h_n\xi_{n-2})_*$, then we see from the proof of Corollary 4.7 (i)
that the only homology class $A \in H_2'$ which has nonzero contributions to
the quantum corrections in (4.4) is $A = A_2$. Thus by (4.4),
$$h \cdot h_p = \cases h_{p+1},       &\text{if $p \le n -2$}\\
                h_n + a_1' \cdot e^{-tn}, &\text{if $p = n -1$}\\
                 h_{n+1} + (a_2'h+a_3'\xi) \cdot e^{-tn}, &\text{if $p = n$}.\\
                \endcases
\eqno (6.2)$$
where $a_1' = \Phi_{(A_2, 0)}(h, h_{n-1}, h_n \xi_{n-1})$,
$a_3' = \Phi_{(A_2, 0)}(h, h_{n}, h_n \xi_{n-2})$, and
$$a_2' = \Phi_{(A_2, 0)}(h, h_{n}, h_{n-1} \xi_{n-1})-c_1 a_3'.$$
By Lemma 3.8, $a_3' = 1$. Thus $a_1 = (a_1'+a_2')$ and the first relation
$f_\omega^1$ in (6.1) is
$$h^{n+1}= ((a_1'+a_2')h + \xi) \cdot e^{-tn} \eqno (6.3)$$
Similarly, from the proof of Corollary 4.7 (ii), we see that
the only homology class $A \in H_2'$ which has nonzero contributions to
the quantum corrections in (4.5) and (4.6) is also $A = A_2$.
By (4.5), $\xi \cdot \xi_p = \xi_{p+1}$
if $p < n -1$, and $\xi \cdot \xi_{n-1} = \xi_{n} + e^{-tn} +
b_2^{(n)} \cdot e^{-tn}$ where
$b_2^{(n)} = \Phi_{(A_2, 0)}(\xi, \xi_{n-1}, h_n \xi_{n-1})$. Thus,
$$\xi^p =
\cases
\xi_{p}, &\text{if $p < n$}\\
\xi_{n} + (1 + b_2^{(n)}) \cdot e^{-tn}, &\text{if $p = n$}\\
\endcases
\eqno (6.4)$$
By (6.2), we have $h \cdot h_p = h_{p+1}$ if $p < n -1$,
and $h \cdot h_{n-1} = h_{n} + b_2^{(0)} \cdot e^{-tn}$ where
$b_2^{(0)} = a_1' = \Phi_{(A_2, 0)}(h_{n-1}, h, h_n \xi_{n-1})$. Thus, we
obtain
$$h^p =
\cases
h_{p}, &\text{if $p < n$}\\
h_{n} + b_2^{(0)} \cdot e^{-tn}, &\text{if $p = n$}\\
\endcases
\eqno (6.5)$$
By (4.6), for $1 \le i \le (n-1)$, $h_{n-i} \cdot \xi_i
= h_{n-i}\xi_i + b_2^{(i)} \cdot e^{-tn}$ where
$b_2^{(i)} = \Phi_{(A_2, 0)}(h_{n-i}, \xi_i, h_n \xi_{n-1})$.
Thus by (6.4) and (6.5), we have
$$h^{n-i} \cdot \xi^i = h_{n-i} \cdot \xi_i =
h_{n-i}\xi_i + b_2^{(i)} \cdot e^{-tn}. \eqno (6.6)$$
Since $\sum_{i=0}^n (-1)^i c_i \cdot h_i \xi_{n-i} = 0$,
it follows from (6.4), (6.5), (6.6) that
$$\sum_{i=0}^n (-1)^i c_i \cdot h^i \cdot \xi^{n-i} =
(1 + \sum_{i=0}^n (-1)^i c_i b_2^{(n-i)}) \cdot e^{-tn}. \eqno (6.7)$$

Next, we compute the above integers $a_1', a_2'$, and $b_2^{(i)}$
where $0 \le i \le n$.

\lemma{6.8} Let $V = T_{\Pee^n}$ with $n \ge 2$ and
$A_2 = (h_{n-1}\xi_{n-1} - nh_{n}\xi_{n-2})_*$.
\roster
\item"{(i)}" $\Phi_{(A_2, 0)}(h, h_n, h_{n-1} \xi_{n-1}) = n$;
\item"{(ii)}" Let $\alpha = h_j\xi_k$ and $\beta = h_s\xi_t$ where
$j, k, s, t$ are nonnegative integers such that
{\rm max}$(j, k) > 0$, {\rm max}$(s, t) > 0$, and $(j+k+s+t) = n$. Then,
$$\Phi_{(A_2, 0)}(\alpha, \beta, h_n \xi_{n-1}) = 1.$$
\endroster
\endproclaim
\noindent
{\it Proof.} (i) By Lemma 2.2 (iv),
$\frak M({A_2}, 0)/PSL(2; \Cee)$ is compact.
By (3.17), we have $h^1(N_f) = 0$ for every $f \in \frak M({A_2}, 0)$.
Thus, $\frak M({A_2}, 0)/PSL(2; \Cee)$ is also smooth.
Fix a line $\ell_0$ in $\Pee^n$. Let $g: \ell_0 \to \Pee(T_{\Pee^n}|_{\ell_0})
\subset \Pee(T_{\Pee^n})$ be the embedding induced by the natural projection
$T_{\Pee^n}|_{\ell_0} = \Cal O_{\ell_0}(1)^{\oplus (n - 1)} \oplus
\Cal O_{\ell_0}(2) \to \Cal O_{\ell_0}(2) \to 0.$
Since $h([g(\ell_0)]) = 1$ and $\xi([g(\ell_0)]) = 2$,
we have $[g(\ell_0)] = [h_{n-1}\xi_{n-1} - (n - 1) h_{n}\xi_{n-2}]_*$.
So $h_{n-1}\xi_{n-1} = [g(\ell_0)]_* + (n - 1) h_{n}\xi_{n-2}$, and
$$\Phi_{({A_2}, 0)}(h, h_{n}, h_{n-1} \xi_{n-1}) =
\Phi_{({A_2}, 0)}(h, h_{n}, [g(\ell_0)]_*)
+ (n - 1) \Phi_{({A_2}, 0)}(h, h_{n}, h_{n}\xi_{n-2}).$$
By Lemma 3.8, it suffices to show that
$\Phi_{(A_2, 0)}(h, h_{n}, [g(\ell_0)]_*) = 1$.
Let $B$ and $C$ be the subvarieties of $\Pee(T_{\Pee^n})$ in general position,
whose homology classes are Poincar\'e dual to $h$ and $h_{n}$ respectively.
Then the homology classes of $\pi(B)$ and $\pi(C)$ in $\Pee^n$ are
Poincar\'e dual to $h$ and $h_n$ respectively.
Let $f: \ell \to \Pee(T_{\Pee^n})$ be a morphism in $\frak M(A_2, 0)$
induced by a surjective map $T_{\Pee^n}|_\ell \to \Cal O_{\ell}(1) \to 0$
for some line $\ell \subset \Pee^n$. If the image $\hbox{Im}(f)$
intersects with $B, C$, and $g(\ell_0)$, then $\ell$ intersects with
$\pi(B)$, $\pi(C)$, and $\pi(g(\ell_0)) = \ell_0$. In other words,
$\ell$ passes through the point $\pi(C)$ and intersects with $\ell_0$.
Moreover, putting $p = \ell \cap \ell_0$ and noticing that
every surjective map $T_{\Pee^n}|_\ell \to \Cal O_{\ell}(1) \to 0$
factors through the natural projection
$T_{\Pee^n}|_\ell = \Cal O_{\ell}(1)^{(n-1)} \oplus \Cal O_{\ell}(2)
\to \Cal O_{\ell}(1)^{(n-1)}$, we conclude that the $(n-1)$-dimensional
subspace $(\Cal O_{\ell}(1)^{(n-1)})|_p$ in
$(T_{\Pee^n}|_\ell)|_p = T_{p, \Pee^n}$ must contain
the $1$-dimensional subspace $(\Cal O_{\ell_0}(2))|_p$ in
$(T_{\Pee^n}|_{\ell_0})|_p = T_{p, \Pee^n}$.
Conversely, let $p \in \ell_0$ and let $\ell_p$ be the unique line
connecting the two points $\pi(C)$ and $p$. If the $(n-1)$-dimensional
subspace $(\Cal O_{\ell_p}(1)^{(n-1)})|_p$ in
$(T_{\Pee^n}|_{\ell_p})|_p = T_{p, \Pee^n}$ contains
the $1$-dimensional subspace $(\Cal O_{\ell_0}(2))|_p$ in
$(T_{\Pee^n}|_{\ell_0})|_p = T_{p, \Pee^n}$,
then there exists a unique surjective map
$T_{\Pee^n}|_{\ell_p} \to \Cal O_{\ell_p}(1) \to 0$
such that the image of the induced morphism $f: \ell_p \to \Pee(T_{\Pee^n})$
intersects $g(\ell_0)$ at the point $g(p)$.
Since there exists a unique point $p \in \ell_0$ such that
the $(n-1)$-dimensional subspace $(\Cal O_{\ell_p}(1)^{(n-1)})|_p$ in
$(T_{\Pee^n}|_{\ell_p})|_p = T_{p, \Pee^n}$ contains
the $1$-dimensional subspace $(\Cal O_{\ell_0}(2))|_p$ in
$(T_{\Pee^n}|_{\ell_0})|_p = T_{p, \Pee^n}$, it follows that
$$\Phi_{(A_2, 0)}(h, h_{n}, [g(\ell_0)]_*) = 1.$$

(ii) It is well-known (see p.176 of \cite{7}) that there is an exact sequence
$$0 \to \Cal O_{\Pee^n} \to \Cal O_{\Pee^n}(1)^{\oplus (n+1)}
\to T_{\Pee^n} \to 0. \eqno (6.9)$$
The surjective map $\Cal O_{\Pee^n}(1)^{\oplus (n+1)} \to T_{\Pee^n} \to 0$
induces the inclusion $\phi: \Pee(T_{\Pee^n}) \subset \Pee^n \times \Pee^n$
such that $\xi$ is the restriction of the $(1, 1)$ class in
$\Pee^n \times \Pee^n$. Let $B, C, q_0$ be the subvarieties of
$\Pee(T_{\Pee^n})$ in general position, whose homology classes are
Poincar\'e dual to $\alpha, \beta, h_n \xi_{n-1}$ respectively.
Then $q_0$ is a point. Put $p_0 = \pi(q_0) \in \Pee^n$.
Now the morphisms in $\frak M(A_2, 0)$ are of the forms
$f: \ell \to \Pee(T_{\Pee^n})$ induced by surjective maps
$T_{\Pee^n}|_\ell \to \Cal O_{\ell}(1) \to 0$
for lines $\ell \subset \Pee^n$. If the image $\hbox{Im}(f)$ passes $q_0$,
then the line $\ell$ passes $p_0$ and $q_0$ is contained in the hyperplane
$$\Pee^{n-2} = \Pee((\Cal O_\ell(1)^{\oplus (n-1)})|_{p_0})
\subset \Pee((T_{\Pee^n}|_\ell)|_{p_0}) = \pi^{-1}(p_0) = \Pee^{n-1}.$$
Conversely, if $\ell$ passes $p_0$ and $q_0$ is contained in the hyperplane
$$\Pee^{n-2} = \Pee((\Cal O_\ell(1)^{\oplus (n-1)})|_{p_0})
\subset \Pee((T_{\Pee^n}|_\ell)|_{p_0}) = \pi^{-1}(p_0) = \Pee^{n-1},
\eqno (6.10)$$
then there exists a unique $f \in \frak M(A_2, 0)$ of the form
$f: \ell \to \Pee(T_{\Pee^n})$ such that $\hbox{Im}(f)$ passes $q_0$;
moreover, putting $q_0 = (p_0, p_0') \in \Pee^n \times \Pee^n$ such that
$\pi$ is the first projection of $\Pee^n \times \Pee^n$,
then $\hbox{Im}(f) = \ell \times \{ p_0' \} \subset \Pee^n \times \Pee^n$.
The set of all lines $\ell$ passing $p_0$ such that $q_0$ is contained
in the hyperplane (6.10) is parameterized by an $(n - 2)$-dimensional
linear subspace $\Pee^{n - 2}$ in $\Pee^n$ (the first factor in
$\Pee^n \times \Pee^n$). It follows that the images
$\hbox{Im}(f) \subset \Pee(T_{\Pee^n})$ sweep a hyperplane
$$H \overset \hbox{def} \to = \Pee^{n - 1} \times \{ p_0' \} \subset
\Pee^n \times \{ p_0' \}. \eqno (6.11)$$
Since $\xi$ is the restriction of the $(1, 1)$ class in $\Pee^n \times \Pee^n$,
$\xi|_H$ is the hyperplane class $\tilde h$ in
$H = \Pee^{n - 1} \times \{ p_0' \} \cong \Pee^{n - 1}$.
Thus $\alpha|_H = \tilde h_{j + k}$ and $\beta|_H = \tilde h_{s+t}$.
Since $(j+k+s+t) = n$ and $B$ and $C$ are in general position,
there is a unique line in $H$ passing $q_0 = (p_0, p_0')$ and
intersecting with $B$ and $C$. Therefore,
$$\Phi_{(A_2, 0)}(\alpha, \beta, h_n \xi_{n-1}) = 1. \qed$$

Finally, we summarize the above computations and prove the following.

\proposition{6.12} The quantum cohomology ring
$H^*_{\omega}(\Pee(T_{\Pee^n}); \Zee)$ with $n \ge 2$ is the ring
generated by $h$ and $\xi$ with the two relations:
$$h^{n+1} = \xi \cdot e^{-tn} \qquad \text{and}
\qquad \sum_{i=0}^n (-1)^i c_i \cdot h^i \cdot \xi^{n-i}
= (1 + (-1)^n) \cdot e^{-tn}.$$
\endproclaim
\proof
By Lemma 6.8 (ii), $a_1' = 1$. By Lemma 3.8, $a_3' = 1$.
By Lemma 6.8 (i),
$$a_2' = \Phi_{(A_2, 0)}(h, h_{n}, h_{n-1} \xi_{n-1}) -c_1 a_3'= -1.$$
Thus by (6.3), the first relation $f_\omega^1$ is
$h^{n+1}= \xi \cdot e^{-tn}$. By Lemma 6.8 (ii),
$b_2^{(i)} = 1$ for $0 \le i \le n$. By (6.7),
the second relation $f_\omega^2$ is
$\sum_{i=0}^n (-1)^i c_i \cdot h^i \cdot \xi^{n-i} =
(1 + \sum_{i=0}^n (-1)^i c_i) \cdot e^{-tn}$. From the exact sequence (6.9),
$c_i = {n+1 \choose i}$ for $0 \le i \le n$.
Therefore, the relation $f_\omega^2$ is
$\sum_{i=0}^n (-1)^i c_i \cdot h^i \cdot \xi^{n-i}
= (1 + (-1)^n) \cdot e^{-tn}$.
\endproof

\enddocument